\documentclass{amsart} \usepackage{comment} \usepackage{amssymb}
\usepackage{enumerate} \usepackage{graphicx}
\newcommand{\unode}[1]{$\stackrel{#1}{\circ}$}
\newcommand{\lnode}[1]{$\stackrel{\textstyle \circ}{\scriptstyle #1}$}

 \def\Simp{{\rm Simp}} \def\rk{{\rm rk}}
\def\tF{\widetilde F} \def\he{\hat{e}} \def\back{{\rm b}}
\def\aback{a^{\back}}  \def\op{{\rm op}}   \def\np{\bigskip\noindent} \def\nl{\smallskip\noindent}
\def\suppress#1{} \def\CoxDiag{M} \def\df{{\small\rm df}} \def\ih{{\small\rm
ih}} \def\rels{{\small\rm rl}} \def\lijntje{\vrule height2.4pt depth-2pt
width0.5in} 
   
\def\vlijntje{\vrule height0.45in depth0.4pt width0.4pt}
\def\vlijn{\buildrel {\hbox to
0pt{\hss$\textstyle\circ$\hss}}\over\vlijntje}  
\def\vtriple#1\over#2\over#3{\mathrel{\mathop{\kern0pt #2}\limits_{\hbox to
0pt{\hss$#1$\hss}}^{\hbox to 0pt{\hss$#3$\hss}}}}
\def\rvtriple#1\over#2\over#3{\mathrel{\mathop{\kern0pt #2}\limits_{\hbox to
0pt{\hss$#3$\hss}}^{\hbox to 0pt{\hss$#1$\hss}}}}

\def\Dnrev{\vtriple{\scriptstyle n}\over\circ\over{}\kern-1pt\lijntje\kern-1pt
\vtriple{\scriptstyle{n-1}}\over\circ\over{}
\cdots\cdots\vtriple{\scriptstyle 4}\over\circ\over{}\kern-1pt\lijntje\kern-1pt
\vtriple{\scriptstyle 3}\over\circ\over{\buildrel
{\scriptstyle 2}\over\vlijn}\kern-1pt\lijntje\kern-1pt
\vtriple{\scriptstyle 1}\over\circ\over{}\kern-1pt}

\def\Esevensigma{\begin{picture}(200,80)
\put(5,4){\lnode{6}}
\put(9,12){\line(1,0){35}}
\put(43,4){\lnode{5}}
\put(47,12){\line(1,0){35}}
\put(81,4){\lnode{4}}
\put(85,12){\line(1,0){35}}
\put(119,4){\lnode{3}}
\put(123,12){\line(1,0){35}}
\put(157,4){\lnode{1}}
\put(161,12){\line(1,0){35}}
\put(195,4){\lnode{0}}
\put(83,14){\line(0,1){35}}
\put(81,48){\unode{2}}
\end{picture}
}

\def\AO{{\mathcal  A}}
\def\MY{{\mathcal  Y}}

\def\A{{\rm A}}
\def\ddA{{\rm A}}
\def\ddD{{\rm D}}
\def\ddE{{\rm E}}

\def\TL{{\rm TL}}

\def\BMW{{\rm B}}
\def\Br{{\rm Br}}

\def\BrM{{\rm BrM}}

\def\D{{\rm D}}
\def\E{{\rm E}}

\def\np{\medskip}
\def\nl{\smallskip\noindent}

\newtheorem{Thm}{\bf{Theorem}}[section]
\newtheorem{Def}[Thm]{\bf{Definition}}
\newtheorem{Defs}[Thm]{\bf{Definitions}}

\newtheorem{Lm}[Thm]{\bf{Lemma}}
\newtheorem{Step}[Thm]{\bf{Step}}

\newtheorem{Prop}[Thm]{\bf{Proposition}}
\newtheorem{Cor}[Thm]{\bf{Corollary}}
\newtheorem{Remark}[Thm]{\bf{Remark}}
\newtheorem{Proof}[Thm]{\em{\hskip-.2mm}}

\newtheorem{Notation}[Thm]{\bf{Notation}}

\renewcommand{\phi}{\varphi}

\DeclareMathOperator{\het}{ht}

\def\Z{{\mathbb Z}}
\def\Q{{\mathbb Q}}

\def\B{{\mathcal B}}
\def\alg{{A}}
\def\He{{\mathcal H}}

\newcommand{\eps}{\varepsilon}

\def\a{\alpha}
\def\b{\beta}
\def\c{\gamma}

\def\homog{\rightsquigarrow}
\def\isog{\leftrightsquigarrow}

\def\cl#1{{#1}^{\rm cl}}

\author{Arjeh M. Cohen
\& David B. Wales}
\address{Arjeh M. Cohen\\
Department of Mathematics and Computer Science\\
Eindhoven University of Technology\\
POBox 513\\
5600 MB Eindhoven\\
The Netherlands}
\email{A.M.Cohen@tue.nl}
\address{David B. Wales\\
Mathematics Department\\
Sloan Lab\\
Caltech\\
Pasadena, CA 91125\\
USA}
\email{dbw@its.caltech.edu}

\title{The Birman--Murakami--Wenzl Algebras of Type $\E_n$}

\date{\today}

\begin{document}

\begin{abstract}
The Birman--Murakami--Wenzl algebras (BMW algebras) of type $\E_n$ for
$n=6,7,8$ are shown to be semisimple and free over the integral domain
$\Z[\delta^{\pm1},l^{\pm1},m]/(m(1-\delta)-(l-l^{-1}))$ of ranks $1,440,585$;
$139,613,625$; and $53,328,069,225$.  We also show they are cellular over
suitable rings.  The Brauer algebra of type $\E_n$ is a homomorphic ring image
and is also semisimple and free of the same rank as an algebra over the ring
$\Z[\delta^{\pm1}]$. A rewrite system for the Brauer algebra is used in
bounding the rank of the BMW algebra above.  The generalized Temperley--Lieb
algebra of type $\E_n$ turns out to be a subalgebra of the BMW algebra of the
same type.  So, the BMW algebras of type $\E_n$ share many structural
properties with the classical ones (of type $\A_n$) and those of type
$\D_n$.
\end{abstract}

\maketitle

\medskip
{\sc keywords:} associative algebra, Birman--Murakami--Wenzl algebra, BMW
algebra, Brauer algebra, cellular algebra, Coxeter group, generalized
Temperley--Lieb algebra, root system, semisimple algebra, word problem in
semigroups

\medskip
{\sc AMS 2000 Mathematics Subject Classification:}
16K20, 17Bxx, 20F05, 20F36, 20M05

\begin{section}{Introduction}\label{intro}
In the paper \cite{CGW} joint with Gijsbers, we introduced
Birman--Murakami--Wenzl algebras (BMW algebras, for short) of simply laced
type, interpreting the classical BMW algebras (introduced in
\cite{BirWen,Mur}) as those of type $\A_n$. Because of the subsequent paper
\cite{CFW}, joint with Frenk, and computations in \cite [Section~$7$]{CGW}
it was expected that these algebras are free of the same rank as the
corresponding Brauer algebras.
This is known for the classical case; see \cite{MorWas}.  In \cite{CGW3}, it
was derived for type $\D_n$.  In this paper, we prove it for types $\E_6$,
$\E_7$, $\E_8$, so that it is established for all spherical simply laced
types.  It is also shown that the algebras are cellular except possibly for
bad primes which are: none for $\A_n$, $2$ for each remaining type, $3$ for
types $\E_n$ $(n=6,7,8)$, and $5$ for $\E_8$.

The classical BMW algebras have a topological interpretation as tangle
algebras; see \cite{MorWas}. In \cite{CGWtangle}, a similar interpretation
was given to BMW algebras of type $\D_n$. Although, in this paper, we
provide bases of the BMW algebras of type $\E_n$ $(n=6,7,8)$ that are built
up from ingredients of the corresponding root systems in the same way as the
other types, an interpretation in terms of tangles is still open.

We use the coefficient ring
$$R = \Z[\delta,\delta^{-1},l,l^{-1},m]/
\left(m(1-\delta)-(l-l^{-1})\right)$$
and recall that, for any simply laced Coxeter diagram $M$,
the BMW algebra $\BMW(M)$ of type
$M$ is the algebra over $R$ given by generators
$g_1,\ldots,g_n,e_1,\ldots,e_n$ and relations as indicated in
Table \ref{BMWTable}.
Here, the indices $i$, $j$, $k$ are nodes of
the diagram $M$. By $i\sim j$ we mean that $i$ and $j$ are adjacent in $M$,
and by $i\not\sim j$ that they are non-adjacent (including the possibility
that they are equal).

\begin{center}
\begin{table}[ht]
\begin{tabular}{|lcl|}
\hline
&& for $i$\\
\hline
(RSrr)&\quad&$g_i^2=1-m(g_i-l^{-1}e_i)$\\
(RSer)&\quad&$e_ig_i=l^{-1}e_i$\\
(RSre)&\quad&$g_ie_i=l^{-1}e_i$\\
(HSee)&\quad&$e_i^2=\delta e_i$\\
\hline\hline
&&for $i\not\sim j$\\
\hline
(HCrr)&\quad&$g_ig_j=g_jg_i$\\
(HCer)&\quad&$e_ig_j= g_je_i$\\
(HCee)&\quad&$e_ie_j= e_je_i$\\
\hline\hline
&&for $i\sim j$\\
\hline
(HNrrr)&\quad&$g_ig_jg_i=g_jg_ig_j$\\
(HNrer)&\quad&$g_je_ig_j=g_ie_jg_i+ m(e_jg_i - e_ig_j + g_ie_j - g_je_i)
 + m^2(e_j-e_i)$\\
(RNrre)&\quad&$g_jg_ie_j=e_ie_j$\\
(RNerr)&\quad&$e_ig_jg_i=e_ie_j$\\
(HNree)&\quad&$g_je_ie_j=g_ie_j+m(e_j-e_ie_j)$\\
(RNere)&\quad&$e_ig_je_i=le_i$\\
(HNeer)&\quad&$e_je_ig_j=e_jg_i+m(e_j-e_je_i)$\\
(HNeee)&\quad&$e_ie_je_i=e_i$\\
\hline
\end{tabular}
\smallskip
\caption{BMW Relations Table, with $i$ and $j$ nodes of $M$}
\label{BMWTable}
\end{table}
\end{center}

\begin{Thm}\label{th:main}
Let $M$ be a simply laced spherical Coxeter diagram.
\begin{enumerate}[(i)]
\item The BMW algebra $\BMW(M)$ is free
of the same rank as the Brauer algebra of type $M$.
\item When tensored with
$\Q(l,\delta)$,  this algebra is semisimple.
\item When tensored with
an integral domain containing
inverses of all bad primes,
$\BMW(M)$ is cellular.
\end{enumerate}
\end{Thm}

Here, the Brauer algebra of type $M$, denoted $\Br(M)$, is as in
\cite{CFW}.  This means it is the free algebra over $\Z[\delta^{\pm1}]$
generated by $r_1,\ldots,r_n,e_1,\ldots,e_n$, with defining relations as
given in Table \ref{BrauerTable} (with the same conventions for $\sim$ and
$\not\sim$). The classical Brauer algebra on Brauer diagrams having $2(n+1)$
nodes and $n+1$ strands introduced in \cite{Brauer} coincides with
$\Br(\A_n)$.  In \cite{CFW} it is shown that $\Br(M)$ is a free
$\Z[\delta^{\pm1}]$-module.  $\Br(M)$ is the image of the ring homomorphism
$\mu: \BMW(M)\to \Br(M)$ sending $e_i$ to $e_i$, and $g_i$ to $r_i$, whilst
specializing $l$ to $1$ and $m$ to $0$.

The ranks $\rk(\Br(M))$ are given in \cite[Table 2]{CFW}; these are
$1,440,585$ for $M =\E_6$, $139,613,625$ for $M =\E_7$, and $53,328,069,225$
for $M =\E_8$, respectively.  Particularly nice bases are provided, which
are parameterized by triples $(B,h,B')$ where $B$ and $B'$ are in the same
orbit $Y$ of special (the technical word being {\em admissible\/}) sets of
mutually orthogonal roots under the Coxeter group $W(M)$ of type $M$ and $h$
belongs to the Coxeter group $W(M_Y)$ whose type $M_Y$ depends only on
$Y$. In the familiar case $M=\A_{n-1}$, the usual basis
consists of Brauer diagrams having $n$ strands; the sets $B$ and $B'$
determine the top and bottom of the Brauer diagram on $n$ strands, where top
and bottom mean the collections of horizontal strands between nodes at the
top and bottom, respectively, and $h$ determines the permutation
corresponding to the vertical strands on the remaining part of the Brauer
diagram (elements of the Coxeter group of type $M_Y = \A_{n-2|B|-1}$).

The generators $e_1,\ldots,e_n$, together with the identity, of the BMW
algebra $\BMW(M)$ satisfy the relations of the Temperley--Lieb algebra of
type $M$ as introduced in Graham's PhD thesis \cite{Gra}.  These are just
the relations (HSee), (HCee), and (HNeee) of Table \ref{BMWTable}.
Therefore $e_1,\ldots,e_n$ together with the identity generate a subalgebra
of $\BMW(M)$ that is a homomorphic image of the Temperley--Lieb algebra over
$R$.  In fact it is the Temperley--Lieb algebra:

\begin{Prop}\label{templieb}
Let $M$ be a simply laced spherical Coxeter diagram. The subalgebra of
$\BMW(M)$ generated by $e_1,\ldots,e_n$ together with the identity is
isomorphic to the Temperley--Lieb algebra of type $M$ over $R$.
\end{Prop}

In particular, the restriction of the ring homomorphism $\mu$ to the
subalgebra of $\BMW(M)$ generated by $e_1,\ldots,e_n$ preserves ranks and maps
a copy of the Temperley--Lieb algebra over $R$ to a copy over
$\Z[\delta^{\pm1}]$.

As mentioned for Theorem \ref{th:main}, this theorem and Proposition
\ref{templieb} are known for $M=\A_n$ (see \cite{MorWas}) and for $M=\D_n$
(see \cite{CGW3}).  The results follow immediately from the results for
connected diagrams $M$ so here only $M=\E_n$ ($n=6,7,8$) need be considered.
The proof of Proposition \ref{templieb} for $M=\E_n$ is given in
\ref{proof1.2}. It rests on the irreducible representations of the
Temperley--Lieb algebras determined by Fan in \cite{Fan}.  Our proof of
Theorem \ref{th:main}(i) for $M=\E_n$ uses Proposition \ref{templieb} as a
base case.  It also uses the special case of \cite[Proposition 4.3]{CGW3}
formulated in Proposition \ref{prop:BMWbasis} below and the rewriting result
stated in Theorem \ref{th:nottheta} further below. It makes use of some
computations in GAP \cite{GAP} for verifications that all possible rewrites
have been covered.

The outline of the paper is as follows. All notions needed for the main
results as well as the main technical results needed for their proofs, are
given in Section \ref{sec:details}.  Section \ref{sec:WC} analyses
centralizers of idempotents occurring in Brauer algebras of type $M = \E_n$
$(n=6,7,8)$. Sections \ref{sec:aB} and \ref{sec:aBcharacterization} together
form the major part of our proof of Theorem \ref{th:main}(i). It runs by
induction on objects from the root system of type $M$, whereas the base
case, related to Temperley--Lieb algebras, is treated in \ref{proof1.2} of
Section \ref{sectiontemplieb}.  The completion of the proof of Theorem
\ref{th:main} as well as a concluding remark is given in Section
\ref{sec:conclusion}.

\section{Detailed statements}
\label{sec:details}
In this section, we describe in detail the statements of the previous
section, the rewrite strategy for their proofs, and the structure of the
Brauer monoid.

Throughout this paper, $F$ is the direct product of the free monoid on
$$r_1,\ldots,r_n,e_1,\ldots,e_n$$ and the free group on $\delta$. Furthermore,
$\pi: F\to \Br(M)$ is the homomorphism of monoids sending each element of
the subset $\{r_1,\ldots,r_n,e_1,\ldots,e_n,\delta,\delta^{-1}\}$ of $F$ to
the element with the same name in $\Br(M)$. Similarly, $\rho: F\to \BMW(M)$
is the homomorphism of monoids sending each element of the subset
$\left\{e_1,\ldots,e_n,\delta,\delta^{-1}\right\}$ of $F$ to the element
with the same name in $\BMW(M)$ and each $r_i$ to $g_i$ $(i=1,\ldots,n)$.  It
follows from these definitions that $\pi = \mu\circ \rho$.

\begin{Defs}\label{df:reduction} \rm
Elements of $F$ are called {\em words}.  A word $ a\in F$ is said to be of
{\em height} $t$ if the number of $r_i$ occurring in it is equal to $t$; we
denote this number $t$ by $\het( a)$.  We say that $ a$ is {\em reducible} to
another word $ b$, that $a$ {\em can be reduced} to $b$, or that $ b$ {\em is
a reduction} of $a$, if $ b$ can be obtained by a sequence of specified
rewrites, listed in Table \ref{BrauerTable}, starting from $a$, that do not
increase the height.  We call a word in $F$ {\em reduced} if it cannot be
further reduced to a word of smaller height.  Following \cite{CGW3}, we have
labelled the relations in Table \ref{BrauerTable} with R or H according to
whether the rewrite from left to right strictly lowers the height or not
(observe that the height of the right hand side is always less than or equal
to the height of the left hand side).  If the number stays the same, we call
it H for homogeneous. Our rewrite system will be the set of all rewrites in
Table \ref{BrauerTable} from left to right and vice versa in the homogeneous
case and from left to right in case an R occurs in its label.  We write
$a\homog b$ if $a$ can be reduced to $b$; for example (RNere) gives
$e_1e_2r_3e_2\homog e_1e_2$ if $2\sim 3$.  If the height does not decrease
during a reduction, we also use the term {\em homogeneous reduction} and write
$a \isog b$; for example, (HNeee) gives $e_2r_1\isog e_2e_3e_2r_1$ if $2\sim
3$.
\end{Defs}

\begin{center}
\begin{table}[ht]
\begin{tabular}{|lcl|lcl|}
\hline
label&\quad&relation&label&\quad&relation\\
\hline
(H$\delta$)&\quad&$\delta$ is central&(H$\delta^{-1}$)&\quad&$\delta\delta^{-1}=1$\\
\hline
&\multispan{4}{\hfill for $i$\hfill}&\\
\hline
(RSrr)&\quad&$r_i^2=1$&
(RSer)&\quad&$e_ir_i=e_i$\\
(RSre)&\quad&$r_ie_i=e_i$&
(HSee)&\quad&$e_i^2=\delta e_i$\\
\hline
&\multispan{4}{\hfill for $i\not\sim j$\hfill}&\\
\hline
(HCrr)&\quad&$r_ir_j=r_jr_i$&
(HCer)&\quad&$e_ir_j= r_je_i$\\
(HCee)&\quad&$e_ie_j= e_je_i$&&&\\
\hline
&\multispan{4}{\hfill for $i\sim j$\hfill}&\\
\hline
(HNrrr)&\quad&$r_ir_jr_i=r_jr_ir_j$&
(HNrer)&\quad&$r_je_ir_j=r_ie_jr_i$\\
(RNrre)&\quad&$r_jr_ie_j=e_ie_j$&
(RNerr)&\quad&$e_ir_jr_i=e_ie_j$\\
(HNree)&\quad&$r_je_ie_j=r_ie_j$&
(RNere)&\quad&$e_ir_je_i=e_i$\\
(HNeer)&\quad&$e_je_ir_j=e_jr_i$&
(HNeee)&\quad&$e_ie_je_i=e_i$\\
\hline
&\multispan{4}{\hfill for $i\sim j\sim k$\hfill}&\\
\hline
(HTeere)&\quad&$e_je_ir_ke_j=e_jr_ie_ke_j$&
(RTerre)&\quad&$e_jr_ir_ke_j=e_je_ie_ke_j$\\
\hline
\end{tabular}
\smallskip
\caption{Brauer Relations Table, with $i$, $j$, and $k$ nodes of $M$}
\label{BrauerTable}
\end{table}
\end{center}

\begin{Prop}\label{prop:BMWbasis}
Let $M$ be of type $\E_n$ for $n\in\{6,7,8\}$.  Let $T$ be a set of words in
$F$ whose image under $\pi$ is a basis of $\Br(M)$.  If each word in $F$
can be reduced to a product of an element of $T$ by a power of $\delta$,
then $\rho(T)$ is a basis of $\BMW(M)$.
\end{Prop}

This proposition is a special case of \cite[Proposition~4.3]{CGW3}.  In view
of this result, Theorem \ref{th:main}(i) follows from Theorem
\ref{th:UniqueReduction} below, which is a rewriting result on the Brauer
monoid \BrM(M) in which computations are much easier than in the
corresponding BMW algebra.
Here, we recall from \cite{CFW}, the Brauer monoid $\BrM(M)$ is the
submonoid generated by $\delta,\delta^{-1}, r_1,\ldots,r_n,e_1,\ldots,e_n$
of the multiplicative monoid underlying the Brauer algebra $\Br(M)$.

Homogeneous reduction, $\isog$, is an equivalence relation, and even a
congruence, on $F$, to which we will refer as {\em homogeneous equivalence}.
We denote the set of its equivalence classes by $\tF$.  Note that
concatenation on $F$ induces a well-defined monoid structure on $\tF$ and
that reduction on $F$ carries over to reduction on $\tF$.

\begin{Thm}\label{th:UniqueReduction}
For $M$ of type $\E_n$ for $n\in \{6,7,8\}$, each element of $\tF$ reduces
to a unique reduced element.
\end{Thm}

The image of $F$ under the homomorphism $\pi$ coincides with
$\BrM(M)$. As $\pi$ is constant on homogeneous equivalence classes, there is
no harm in interpreting $\pi$ as a map $\tF\to\BrM(M)$.  Let $T_\delta$ be
the set of reduced words in $\tF$. By definition of $\BrM(M)$ and Theorem
\ref{th:UniqueReduction}, the restriction of $\pi$ to $T_\delta$ is a
bijection onto $\BrM(M)$.  The cyclic group generated by $\delta$ acts
freely by multiplication on $T_\delta$.  Choose $T$ to be a set of
representatives in $T_\delta$ for this action.  As $\pi$ is equivariant with
respect to this action and $\Br(M)$ is canonically isomorphic to the free
$\Z$-algebra over $\BrM(M)$, the restriction of $\pi$ to $T$ is a bijection
onto a basis of $\Br(M)$ over $\Z[\delta^{\pm1}]$.  Consequently,
Proposition \ref{prop:BMWbasis} applies, giving that $\rho(T)$ is a basis of
$\BMW(M)$. This reduces the proof of Theorem \ref{th:main}(i) to a proof of
Theorem \ref{th:UniqueReduction}. We shall however prove a stronger version
of the latter theorem in the guise of Theorem \ref{th:nottheta}.

\medskip
We next describe the set $T_\delta$ of reduced words in $\tF$.  Our starting
point is a finite set, denoted $\AO$ and introduced in
\cite[Section~$3$]{CGW2}, on which the Brauer monoid $\BrM(M)$ acts from the
left.  Elements of $\AO$ are particular, so-called {\em admissible}, sets of
mutually orthogonal positive roots from the root system $\Phi$ of type $M$
(see below for the precise definition).  A special element of $\AO$ will be
the empty set $\emptyset$. By restriction, the Coxeter group $W$ of type $M$
also acts on $\AO$ and we will use a special set $\MY$ of $W$-orbit
representatives in $\AO$, whose members we can associate with subsets $Y$ of
the nodes of $M$ on which the empty graph is induced; such sets of nodes are
called {\em cocliques} of $M$.  The empty coclique of $M$ represents the
member of $\AO$ equal to $\emptyset$, which is fixed by $W$.

Let $Y$ be a coclique of $M$. The element $e_Y$ of $\tF$ denotes the
product over all $i\in Y$ of $e_i$.  As no two nodes in $Y$ are adjacent,
(HCee) implies that the $e_i$ $(i\in Y)$ commute, so it does not matter in
which order the product is taken.  For each node $i$ of $M$, put $\he_i =
e_i\delta^{-1}$ and put $\he_Y = e_Y\delta^{-|Y|} = \prod_{i\in
Y}\he_i$. These are idempotents.

Corresponding to $Y$, there is a unique smallest admissible element of $\AO$
containing $\{\a_i\mid i\in Y\}$, denoted $B_Y$.  With considerable effort,
we are able to define, for each $B$ in the $W$-orbit $WB_Y$ of $B_Y$, an
element $a_B$ of $\tF$ that is uniquely determined up to powers of $\delta$
by $\pi(a_B)\emptyset =\pi(a_B)B_Y = B$ and certain minimality conditions.  The
precise statements appear in Theorem \ref{prop:aB} below. Also, we will
identify a subset $T_Y$ of $\tF$ of elements commuting with $e_Y$ in $\tF$
and in bijective correspondence with a Coxeter group of type $M_Y$; see
Proposition \ref{th:WC} and Table \ref{table:types2}. Now
\begin{eqnarray} \label{dfT}
T_\delta &=& \Big\{ \delta^i a_B \he_Y h a_{B'}^{\op}\Big.\,\Big|
\,  Y\in {\MY};\, B,B'\in WB_Y;\, h\in
T_Y, i\in\Z
\Big\}.
\end{eqnarray}
Here the map $a\mapsto a^{\op}$ on $F$ is obtained (as in \cite[Notation
3.1]{CGW3}) by replacing an expression for $a$ as a product of its
generators by its reverse. This induces an antiautomorphism on $\tF$ and
on $\BrM(M)$.  Equality (\ref{dfT}) illustrates how the triples $(B,h,B')$
alluded to before parameterize the elements of $T$. The detailed
description of $T$ reveals a combinatorial structure that will be used to
prove the semisimplicity and cellularity parts of Theorem
\ref{th:main} (see Section \ref{sec:conclusion}).

\np
We now give precise definitions of the symbols introduced for the
description of $T$. Throughout this section, we let $M$ be a connected
simply laced spherical diagram.  Instead of $W(M)$ we also write $W$ for the
Coxeter group of type $M$.

The combinatorial properties of the root system $\Phi$ of type $M$ that we
will discuss here are crucial.  We first recall the definition of
admissible. A set $X$ of orthogonal positive roots is called {\em
admissible} if, for any positive root $\b$ of $\Phi$ that has inner
product~$\pm1$ with three roots, say $\b_1$, $\b_2$, $\b_3$, of $X$, the sum
$2\b-\sum_{i=1}^3(\b,\b_i)\b_i$ is also in $X$.  In \cite{CFW} and
\cite{CGW2} it is shown that any set $X$ of orthogonal positive roots is
contained in a unique smallest admissible set, which is called its {\em
admissible closure} and denoted $\cl{X}$.  Now $W$ acts elementwise on
admissible sets with the understanding that negative roots are being
replaced by their negatives: for $w\in W$ and $B\in\AO$, we have $wB = \{
\pm w\a\mid\a\in B\}\cap\Phi^+$.  If $M = \A_n$, all sets of mutually
orthogonal positive roots are admissible.

In \cite{CGW2}, a partial ordering $<$ with a single maximal element is
defined for each $W$-orbit in $\AO$. An important property of this partial
ordering is that, if $i$ is a node of $M$ and $B\in\AO$, then $r_iB<B$ is
equivalent to the existence of a root $\b$ of minimal height in $B\setminus
r_iB$ for which $\het(r_i\b)<\het(\b)$; see \cite[Section~$3$]{CGW2}. A
useful property of this ordering is that, for each $i$ and $B$, the sets $B$
and $r_iB$ are comparable.  The definition of $M_Y$ depends on this
ordering.  The ordering is also involved in a notion of height for elements
of $\AO$, denoted $\het(B)$ for $B\in\AO$, which satisfies $\het(B)<\het(C)$
whenever $B,C\in\AO$ satisfy $B<C$. Moreover, if $r_iB>B$, then $\het(r_iB)
= \het(B)+1$. (See Definitions \ref{eY} below for further details.)


\medskip
Nonempty representatives of $W$-orbits in $\AO$ are listed in
\cite[Table~$2$]{CGW2} and, for $M=\E_n$ $(n=6,7,8)$, in
Table~\ref{table:types2}.  Each line of Table~\ref{table:types2} below the
header corresponds to a single $W$-orbit in $\AO$.

\begin{Defs}\label{df:HY}
\rm
By $\MY$ we denote the set consisting of
the empty set and the cocliques $Y$ of $M$ listed in column 5 of Table
\ref{table:types2}.

Let $Y\in \MY$.  We recall that $B_Y=\cl{\{\a_i\mid i \in Y\}}$, the
admissible closure of the set of simple roots indexed by $Y$.  It is a fixed
representative of a $W$-orbit in $\AO$.  The Coxeter type $M_Y$ is the
diagram induced on the nodes of $M$ whose corresponding roots are orthogonal
to all members of the single maximal element of
$WB_Y$ with respect to the partial order $<$ (see \cite{CGW2}, where the type is denoted
$C_{WB_Y}$).

We denote by $H_Y$ the subsemigroup of $\tF$ generated by the elements of
$S_Y$ and $\he_Y$ occurring in the sixth column of Table~\ref{table:types2}.
Finally, we write $T_Y$ for the subset of $\tF$ consisting of reduced
elements of $H_Y$.
\end{Defs}

We will show that $H_Y$ is a monoid with identity $\he_Y$ whose generators
$S_Y$ satisfy certain Coxeter relations.  Then $\pi$ maps $H_Y$ onto a
quotient of the Coxeter group of type $M_{Y}$. In fact, in Proposition
\ref{th:WC} the image $\pi(H_Y)$ turns out be isomorphic to the Coxeter
group, and $T_Y$ turns out to be in bijective correspondence with
$W(M_{Y})$.

The first column of Table~\ref{table:types2} indicates to which type $M$ the
row belongs. By now the meaning of the fifth column (the coclique $Y$ of
$M$), the second column (the size of $B_Y$), fourth column (the type $M_Y$),
and the one but last column (a distinguished subset $S_Y$ of $\tF$),
should be clear. We describe the other columns of this table.

The third column lists the Coxeter type of the root system on the roots
orthogonal to $B_Y$.  The centralizer $C_W(B_Y)$ of $B_Y$ in $W$ is analyzed
in \cite{CGW2}. It is the semi-direct product of the elementary abelian
group of order $2^{|B_Y|}$ generated by the reflections in $W$ with roots in
$B_Y$ and the subgroup $W(B_Y^\perp\cap\Phi)$ of $W$ generated by
reflections with roots in $B_Y^\perp\cap\Phi$.  The normalizer, or setwise
stabilizer, $N_W(B_Y)$ of $B_Y$ in $W$ can be larger and is described in
\cite[Table~$1$]{CGW2}.



The last column lists the sizes of the collections, $(WB_Y)^0$, of
admissible sets of height~$0$ in the $W$-orbit $WB_Y$ of $B_Y$. This data
will not be needed until Section \ref{sectiontemplieb}.

\begin{table}[ht]
\begin{center}
\begin{tabular}{|cc|ccccc|}
\hline ${\CoxDiag}$ & $|B_Y|$ & $B_Y^\perp$ & $M_{Y}$ &$Y$&
$S_Y=\{x\he_Y\mid x$ as below $\}$&$|(WB_Y)^0|$\\ \hline
$\E_6$ & $1$ & $\A_5$ & $\A_5$ &$6$&$e_6e_5e_4e_3r_2e_4e_5,r_1,r_2,r_3,r_4$&6\\
$\E_6$ & $2$ & $\A_3$ & $\A_2$ &$4,6$&$e_4e_3r_2,r_1$&20\\
$\E_6$ & $4$ & $\emptyset$ &  $\emptyset$&$2,3,6$&- &15\\
\hline
$\E_7$ & $1$ & $\D_6$ & $\D_6$ &$7$&$e_7\cdots e_3r_2e_4e_5e_6,r_1,\ldots,r_5$&7\\
$\E_7$ & $2$ & $\A_1 \D_4$ & $\A_1 \A_3$
&$5,7$&$e_5e_4e_3r_2e_4, r_1$&27\\
$\E_7$ & $3$ & $\D_4$ & $\A_2$ &$2,5,7$&$r_1, r_3$&21\\
$\E_7$ & $4$ & $\A_1^3$ & $\A_1$ &$2,3,7$&$r_5$&35\\
$\E_7$ & $7$ & $\emptyset$ & $\emptyset$ & $2,3,5,7$&-&15\\
\hline
$\E_8$ & $1$ & $\E_7$ & $\E_7$ & $8$&$e_8\cdots e_3r_2e_4\cdots e_7,r_1,\ldots,r_6$&8\\
$\E_8$ & $2$ & $\D_6$ & $\A_5$ & $6,8$&$e_6e_5e_4e_3r_2e_4e_5,r_1,r_2,r_3,r_4$&35\\
$\E_8$ & $4$ & $\D_4$ & $\A_2$ & $2,3,8$&$r_5,r_6$&84\\
$\E_8$ & $8$ & $\emptyset$ & $\emptyset$ &$2,3,5,8$&-&50\\
\hline
\end{tabular}
\end{center}
\smallskip
\caption{Nonempty cocliques $Y$ of $M$ and
admissible sets $B_Y$.}
\label{table:types2}
\end{table}

\medskip
As a result of this description of the reduced element set $T_\delta$ in
(\ref{dfT}), the size of $T_Y$ coincides with $|W(M_{Y})|$ and the rank of
$\Br(M)$ over $\Z[\delta^{\pm1}]$ is
\begin{eqnarray*}
 |T|
&=& \sum_{Y\in {\MY}} |W(M_{Y})|\cdot |WB_Y|^2.
\end{eqnarray*}
Substituting the data of Table \ref{table:types2},
we find the values of \cite[Table 2]{CFW} (and
listed above Proposition \ref{templieb}).
This description is a strengthening of \cite[Proposition~$4.9$]{CFW}.

We continue by recalling the action of the monoid $\Br(M)$ on $\AO$
introduced in \cite{CFW}.

\begin{Def}\label{def:Baction} \rm
Let $M$ be a simply laced spherical Coxeter diagram and let $\AO$ be the union
of all $W$-orbits of admissible sets of orthogonal positive roots (so the
empty set is a member of $\AO$).  The action of $W$ on $\AO$ is as discussed
above. The action of $\delta$ is taken to be trivial, that is $\delta (X) = X$
for $X\in \AO$.  This action extends to an action of the full Brauer monoid
$\BrM(M)$ determined as follows on the remaining generators, where $i$ is a
node of $\CoxDiag$ and $B\in \AO$.
\begin{equation}\label{def:sigma_e}
e_iB = \begin{cases}B&\mbox{ if  } \a_i\in B,\\
\cl{(B\cup\{\a_i\})}&\mbox{ if  }\a_i\perp B,\\
r_{\b}r_i B &\mbox{ if  }\b\in B\setminus\a_i^\perp.
\end{cases}
\end{equation}
It is shown in \cite[Theorem~3.6]{CFW} that this is an action.
\end{Def}

Using the antiautomorphism $a\mapsto a^{\op}$ we obtain a right action of
$\BrM(M)$ on $\AO$ by stipulating $Ba = a^{\op}B$ for $B\in\AO$ and
$a\in\BrM(M)$.  (We will also write $a^{\op}$ for the reverse of a word $a$ in
$F$ or of an element $a$ of $\tF$.)

\begin{Defs}\label{eY}  \rm
As indicated above, by $B_Y$ we denote the admissible closure of $\{\a_i
\mid i\in Y\}$.  It is a minimal element of the poset on $WB_Y$ induced by
the partial ordering $<$ defined on $\AO$. If $d$ is the distance in
the Hasse diagram for $WB_Y$ from $B_Y$ to the unique maximal element of
$WB_Y$ (whose existence is proved in \cite[Corollary 3.6]{CGW2}), then, for
$B\in WB_Y$, the {\em height} of $B$, notation $\het(B)$, is $ d - \ell $,
where $\ell$ is the distance in the Hasse diagram from $B$ to the maximal
element. In particular, $\het(B_Y)=0$ and the maximal element has height
$d$.

The {\em level} of an admissible set $B$, notation $L(B)$, is the pair
consisting of the height of $B$ and the multiset $\{\het(\b)\mid \b\in B\}$.
These are ordered by first height of $B$ and then lexicographically, with the
lower heights of roots of $B$ coming first.

For any given $B\in \AO$ we define
$\Simp(B)$ to be the set of simple roots in $B$.
\end{Defs}

\np Our proof of Theorem \ref{th:UniqueReduction} consists of the following
reduction strategy.  Let $a\in\tF$. Then $B = \pi(a)\emptyset$ and $B' =
\emptyset \pi(a)$ belong to the same $W$-orbit of $\AO$.  Fix $Y\in{\MY}$ be
such that $ B\in W B_Y$.  We will show $a\homog \delta^i a_B \he_Y h
a_{B'}^{\op}$ for some $h\in H_Y$ and $i\in\Z$.  By using the
Matsumoto--Tits rewrite rules for Coxeter groups, cf.~\cite{Matsum,Tits69},
we may even take $h\in T_Y$ (cf.~Definitions
\ref{df:HY}).  In summary, with $T_\delta$ as in (\ref{dfT}),
the proof of Theorem \ref{th:UniqueReduction} is a direct consequence of the
theorem below.  Recall that $T_Y$ is the set of reduced element of $H_Y$.

\begin{Thm}\label{th:nottheta}
Let $M$ be a simply laced spherical Coxeter diagram.  Suppose that $ a$ is a
word in $F$.  Let $Y\in{\MY}$ be such that $B_Y$ and $B = \pi(a)\emptyset$
are in the same $W$-orbit. Then $B'=\emptyset\pi(a)$ is in the same
$W$-orbit as $B$ and $B_Y$, and $a\homog \delta^i a_B \he_Y h a_{B'}^{\op}$
for some $i\in\Z$ and $h\in T_Y$. In particular, each element of $\tF$
reduces to a unique element of $T_\delta$, and each element of $T_\delta$ is
reduced.
\end{Thm}

By \cite[Proposition 4.9]{CFW} and the rank computations in [loc.~cit.], the
monomials $\pi(a_B \he_Y h a_{B'}^{\op})$ in $\Br(M)$ are indeed distinct
for distinct triples $(B,h,B')$, as are their multiples by different powers
of $\delta$. So the burden of proof is in the uniqueness of $a_B$ and $h$
when given $a$ with $B = \pi(a)\emptyset$.  The proof of Theorem
\ref{th:nottheta} is presented in \ref{proof:nottheta} and is based on the
three main results, Theorems \ref{prop:aB}, \ref{th:WC},
\ref{th:aBcharacterization}, which are stated below.

\begin{Cor}\label{cor:nottheta} Under the hypothesis of
Theorem~\ref{th:nottheta}, if $a$ and $a'$ are two words of height
$\het(aB_Y)$ with $aB_Y = a'B_Y$, then $a\isog a'$ up to powers of $\delta$.
\end{Cor}

We now introduce an algorithm that will give, for any given $B\in \AO$,
a word $a_B$ having the required properties for the definition of $T$.  We
also introduce another word $\aback_B$, which moves $B$ to $B_Y$ (as defined
in Theorem \ref{th:nottheta}).  We need certain words, called Brink--Howlett
words, from the subsemigroup of $\tF$ generated by $e_1,\ldots,e_n$ that
are specified in Definition \ref{df:HB}. They originate from
\cite{BrinkHowlett} and were also described for reflection groups in the
earlier paper \cite{How}.  The Brauer elements of these Brink--Howlett words
have the property that, whenever $Y$ and $Y'$ are two cocliques of $M$ with
$|Y|=|Y'|$ such that $B_Y$ and $B_{Y'}$ are in the same $W$-orbit, then they
move one to the other in the $\BrM(M)$-action on $\AO$.

\begin{Def}\label{def:aB}
\rm For $B\in WB_Y$, we denote by $a_B$, respectively $\aback_B$, a word in
$\tF$ constructed according to the following rules.
\begin{enumerate}[(i)]
\item
If $|\Simp(B)|=|\Simp(B_Y)|$, then $a_B$ is the Brink--Howlett word that, in
the left action, takes
$B_Y$ to $B$, followed by $\he_Y$.  Moreover, $\aback_B$ is the
Brink--Howlett word taking $B$ to $B_Y$ in the right action,
followed by $\he_Y$.
\item
If $r_kB<B$ for some node $k$, then $a_{B} = r_ka_{r_kB}$
and  $\aback_{B} = r_k\aback_{r_kB}$.

\item
Otherwise, there are adjacent nodes $j$ and $k$ of $M$ with $\a_j\in B$ such
that
$\het(e_kB) = \het(B)$ and $L(e_kB)<L(B)$.
Then $a_{B} = e_ja_{e_kB}$
and  $\aback_{B} = e_k \aback_{e_kB}$.
\end{enumerate}

The nodes $k$ described in (iii) are called {\em
lowering-e-nodes} for $B$.  The nodes $k$ for which $r_kB<B$ are called {\em
lowering nodes} for $B$.
\end{Def}

Notice that $\pi(a_B)\emptyset = \pi(a_B)B_Y = B$ and $ B\pi(\aback_B) = B_Y
$.  Rule (i) only deals with admissible sets of height $0$. The equality of
heights in (iii) for $e_kB$ and $B$ is a consequence of the other properties,
as will be clear from Lemma \ref{sameheight}.

The only rule changing the height in the poset $\AO$ is (ii) and here it is
lowered by exactly by $1$.  This also means $a_B$ is reduced as each $r_k$
in (ii) lowers the height of $a_B$ as well as the height of $B$ by $1$ so
there must be at least $\het(B)$ occurrences of $r_k$'s in any word $a\in F$
with $\pi(a)\emptyset = B$.  This gives the very important property, stated
in (i) below, relating the heights of $a_B$ and of $B$.

\begin{Prop}\label{property:aB}  For each $B\in\AO$, the following holds.
\begin{enumerate}[(i)]
\item
$\het(B) = \het(a_B)$.
\item The word $a_B$ is reduced.
\item
There exist words $a_B$ and $a_B^{\rm b}$ in $\tF$ constructed as in
Definition \ref{def:aB}.
\end{enumerate}
\end{Prop}

\begin{proof}
Assertion (i) is a direct consequence of the construction of $a_B$ in
Definition \ref{def:aB}. As any word $a\in F$ with
$\pi(a)B_Y = B$ satisfies
$\het(a)\ge \het(B)$, assertion (ii) follows from (i).
So it remains to establish (iii).

To this end, we verify that the conditions of Definition \ref{def:aB}(iii)
are always satisfied so that words $a_B$ and $a_B^{\rm b}$ constructed as in
Definition \ref{def:aB} are guaranteed to exist.  We know there are no nodes
$k$ for which $r_kB<B$.  If there are fewer than $|Y|$ simple roots in $B$,
take one of minimal height, say $\b$, in $B$ that is not simple and a node
$k$ lowering $\{\b\}$.  As $B$ and $r_kB$ are comparable, we must have
$r_kB>B$, and so there is a node $j$ for which $\a_j\in B$ is raised by $k$
and so $k\sim j$.  Now $e_kB=r_jr_kB$ has height $\het(B)$.  Under the
action of $e_k$, the simple root $\a_j$ in $B$ is replaced by the simple
root $\a_k$ in $e_kB$, and $\b$ is replaced by $\b-\a_k-\a_j$, so
$L(e_kB)<L(B)$ unless there is a node $i\sim k$ with $\a_{i}$ also in $B$.
In the latter case we use the fact that $B$ is admissible, which implies
$\b-\a_j-\a_{i}-2\a_k$ also belongs to $B$. As its height is lower than
$\het(\b)$, it must be simple.

So we may assume that $B$ has at least three simple roots. We are done in
the case of sets of size at most $4$.
Admissible sets $B$ of size $7$ or
$8$ in $\E_7$ and $\E_8$ remain.
In these cases, take $\b'$ in $B\setminus \cl{\Simp(B)}$ of minimal height
and take a node $k'$ lowering $\b'$.  Then $k'\sim l$ for at most one node $l$
with $\a_l\in \Simp(B)$.  This $k'$ will be as required.
\end{proof}

\begin{Thm}\label{prop:aB}
Let $M\in \{\E_6,\E_7,\E_8\}$ and $Y\in{\MY}$.
For each $B\in W B_Y$ there is, up to homogeneous
equivalence and powers of $\delta$, a unique word $a_B$ in $F\he_{Y}$
satisfying Definition~\ref{def:aB}. This word has height $\het(B)$ and moves
$\emptyset$ to $B$ in the left action: $\pi(a_B) \emptyset = B$.  Moreover,
there is a word $\aback_B$ in $F$ of height $\het(B)$ that satisfies $B\pi(\aback_B) =
B_Y$.
\end{Thm}

The proof of this result is described after Theorem
\ref{th:aBcharacterization}.  Contrary to $a_B$, the words $\aback_B$ are
not uniquely determined.

If $\het(B) = 0$, then $a_B$ and $\aback_B$ are Temperley-Lieb words as
discussed in Section~\ref{sectiontemplieb}.  Clearly, then $r_kB\ge B$ for
all nodes $k$ of $M$.  The converse is true for $M=\ddA_n$: the word $a_B$
will be a product of an element from $W$ and a Temperley--Lieb word.  For
other types $M$, this is not necessarily the case.  An example is the
admissible set $B =\{\a_4,\a_1+\a_2+2\a_3+2\a_4+\a_5\}$ for $M=\E_6$. As
$r_1$ and $r_4$ leave $B$ invariant and $r_2$, $r_3$, $r_5$, and $r_6$ raise
$B$, there is no lowering node for $B$; consequently $a_B$ cannot begin with
an element from $W$, but its height equals $2$. In fact
we can take $a_B=e_4r_2r_5e_3e_4e_5e_1e_3\he_4\he_6$ and
$\pi(a_B)\emptyset =r_3r_4r_2r_5r_1r_3r_5r_4r_6r_5r_3r_1r_4r_3B_Y$,
with $Y = \{4,6\}$.
In particular, $B$ is an
admissible set as in Case (iii) of Definition \ref{def:aB} with $\het(B)>0$.
In accordance with Proposition \ref{property:aB} the Temperley--Lieb word
$e_3$ satisfies $L(e_3B)<L(B)$ and $e_3B$ has lowering nodes $2$ and $5$.

\begin{Thm}\label{th:WC}
Let $M\in \{\E_6,\E_7,\E_8\}$ and $Y\in\MY$. The Matsumoto--Tits rewrite
rules of type $M_Y$ are satisfied by $S_Y$ in $F$ with respect to $\homog$,
with identity element $\he_Y$.
Moreover, the set $T_Y$ of reduced words of the submonoid $H_Y$ of $\tF$
generated by $S_Y$ are in bijective correspondence with the elements of
$W(M_Y)$.
\end{Thm}

The rewriting for $H_Y$ is handled via the Matsumoto--Tits rewrite rules for
$W(M_{Y})$, the Coxeter group of type $M_{Y}$.  The proof and a further
structure analysis of $H_Y$ is given in \ref{proofWC}.

The rewriting for $a\in\tF$ is handled via the following
behavior of the elements $a_B$ under left multiplication with generators of
$\BrM(M)$.
Observe that $a_B$ ends in $\he_{Y}$.

\begin{Thm}\label{th:aBcharacterization}
Let $M\in \{\E_6,\E_7,\E_8\}$ and $Y\in\MY$.
For each $B\in W B_Y$ the element
$a_B$
of $\tF$ has height $\het(B)$ and satisfies
the following three properties for each node $i$ of $M$.
\begin{enumerate}[(i)]
\item
$r_i a_{B}\homog a_{r_iB}h$ for some $h\in H_Y$.  Furthermore, if $r_iB>B$,
then $h=\he_Y$, the identity in $H_Y$.
\item If $|e_iB| = |B|$, then
$e_ia_B\homog a_{e_iB}h$ for some $h\in \delta^{\Z}H_Y$
and $\het(e_iB)\le \het(B)$.
\item If $|e_iB| > |B|$, then
$e_ia_B$ reduces to an element of $\BrM(M) e_{U}\BrM(M)$
for some set of nodes $U$ strictly containing $ Y$.
\end{enumerate}
\end{Thm}


Fix $M\in \{\E_6,\E_7,\E_8\}$.  The proofs of Theorems~\ref{prop:aB} and
\ref{th:aBcharacterization} are closely related. Actually, the assertions
are proved by induction on the rank of $M$ as well as the level $L(B)$ of
the admissible set $B$ involved.  In Section \ref{sec:aB} we prove the
statement of Theorem~\ref{prop:aB} for $B\in \AO$ assuming the truth of the
statements of both theorems for elements in $\AO$ of level less than $L(B)$.
In Section \ref{sec:aBcharacterization} we prove the statement of
Theorem~\ref{th:aBcharacterization} for $B\in \AO$ assuming the truth of the
statements of Theorem~\ref{prop:aB} for elements in $\AO$ of height less
than or equal to $L(B)$ and of Theorem \ref{th:aBcharacterization} for
elements of height strictly less than $L(B)$. The base case for the
induction, $\het(B) = 0$, is covered by Corollary \ref{cor:aBszero}.
As the results are already
proved for types $\A_n$ and $\D_n$, see \cite[Section 4]{CGW3}, we also
assume the validity of the theorems for BMW algebras whose types have
strictly lower ranks than $M$.

\end{section}

\begin{section}{The Temperley--Lieb Algebra}
\label{sectiontemplieb}

The parts of Theorems~\ref{prop:aB} and \ref{th:aBcharacterization}
concerned with admissible sets $B$ of height zero are proved in this section.
We also provide a proof of Proposition \ref{templieb}.

There are some natural height preserving actions by $e_i$ which arise in many
of our calculations.

\begin{Lm}\label{sameheight}  Let $B\in \B$ and let $j$ be a node of $M$.
Then $\a_j\in e_jB$. Assume further that $i$ is a node of $M$ with $\a_i\in B$
and $i\sim j$. Then  $\het(B)=\het(e_jB)$.  Furthermore,
$B=e_ie_jB$ and $e_jB=e_je_i(e_jB)$.
\end{Lm}

\begin{proof}
The first assertion is direct from the last rule of (\ref{def:sigma_e}) and
the observation that $r_\beta r_j \beta = \a_j$ if $\beta\in B\setminus
\a_j^\perp$.

As for the second assertion, the last rule of (\ref{def:sigma_e}) and
$\a_i\in B\setminus\a_j^\perp$ give $e_jB=r_ir_jB$.  Now $r_jB>B$ as
$r_j\a_i=\a_i+\a_j$, so an element of height $1$ becomes of height $2$.
This means $\het (r_jB)=\het(B)+1$.  No simple root $\a_k\in r_jB$ is raised
to $\a_i+\a_k\in e_jB$, for otherwise we would have $ 0 = (\a_i+\a_k,\a_j) =
-1$, a contradiction.  But $r_i(\a_i+\a_j)=\a_j$ and so an element of height
$2$ in $r_jB$ is lowered to height $1$.  This means $\het(e_jB)=\het(r_jB)-1
= \het(B)$.  As $e_jB$ contains $\a_j$, we find
$e_i(e_jB)=r_jr_ie_jB=r_jr_ir_ir_jB=B$.  Finally, $e_i=e_ie_je_i$ implies
$e_iB=e_ie_j(e_iB)$.
\end{proof}

Each $W$-orbit $\B$ in $\AO$ contains a certain number of admissible sets
$B$ with the maximal number of simple roots, which is $|Y|$ of Table
\ref{table:types2}.  This is the size of $B$ except for sets of size four,
seven, and eight.  For sets of size four, the nodes of these simple roots
can be taken to be $\{2,3,n\}$ and for sets of sizes seven and eight (in
case $\E_7$ as well as $\E_8$) they can be taken to be $\{2,3,5,n\}$.  If
$B$ has the maximal number of simple roots in its $W$-orbit, it is the
admissible closure of $\Simp(B)$.

\begin{Lm}\label{lm:HB}
Let $U$ and $U'$ be two cocliques of $M$ such that
$B_U$ and $B_{U'}$ are in the same $W$-orbit. Then
there is a word $a=e_{i_1}\cdots e_{i_s}$ with
$\pi(a) B_U =  B_{U'}$.
\end{Lm}

\begin{proof}
The work \cite{BrinkHowlett} of Brink--Howlett
shows that $\{\a_i\mid i\in U\}$ can be mapped to
$\{\a_i\mid i\in U'\}$ by a sequence of products $r_{i_t}r_{j_t}$
$(t=1,\ldots,s$) of two reflections with $i_t\sim j_t$
such that $\a_{i_t}\in r_{i_{t-1}}r_{j_{t-1}}\cdots
r_{i_{1}}r_{j_{1}}\{\a_i\mid i\in U\}$.  So $B_{U'} =
r_{i_s}r_{j_s}\cdots r_{i_1}r_{j_1}B_U$ and, for the corresponding intermediate
images $B_{t}=r_{i_t}r_{j_t}\cdots r_{i_1}r_{j_1}B_U$ of $B_U$, the root
$\a_{i_t}$ belongs to $B_{t}$ and $\a_{j_t}$ belongs to $B_{t+1}$, so
$r_{i_t}r_{j_t}B_t$ coincides with $e_{j_t}B_t$. Consequently, the word $a =
e_{j_s}\cdots e_{j_1}$ satisfies $B_{U'} = \pi(a)B_U$, as required.
\end{proof}

\begin{Def}\label{df:HB}
\rm
The words $a$ appearing in Lemma \ref{lm:HB} are called
{\em Brink--Howlett words}.
\end{Def}

\np
These words enter as part of Definition \ref{def:aB} of $a_B$.  The
method is to act by $r_i$ and $e_i$ in such a way as to get the correct
maximum number of simple roots in $B$ and then to act by Brink--Howlett
words to get the fixed one $B_Y$.  In the definition of $a_B$ the action on
the left takes $B_Y$ to $B$.  These other elements of $\B$ are all at height
$0$ by Lemma~\ref{sameheight}.  They are the lowest height possible by the
properties of $a_B$.

\begin{Notation}\label{not:TL}
\rm Let $\TL(M)$ be the subalgebra of $\Br(M)$ generated by the elements
$e_i$ together with the identity in $\Br(M)$. So, by construction it is a
homomorphic image of the Temperley--Lieb algebra of type $M$, that is, the
free algebra with identity generated by $e_i$ $(i=1,\ldots,n)$ subject to
the relations (HSee), (HCee), and (HNeee).
\end{Notation}

In Proposition \ref{prop:irreducible} we prove that $\TL(M)$ is isomorphic
to the Temperley--Lieb algebra of type $M$.  Up to powers of $\delta$, the
monomials in $\TL(M)$ are elements of the form $e_{k_1}\cdots e_{k_l}$.

\begin{Notation}\label{df:height0}
\rm For a given $Y\in\MY$,
we denote the collection of height $0$ sets in $WB_Y$ by $(WB_Y)^0$.
\end{Notation}

\begin{Lm}\label{height0aB}  Let $Y\in\MY$
and $B\in (WB_Y)^0$.  Then $a_B$ is a product $e_{i_1}e_{i_2} \cdots
 e_{i_r}\he_Y$ such that each $i_{j+1}$ is adjacent to a node
associated with a simple root in $e_{i_j}\cdots e_{i_r}B_Y$.
Also $\aback_B$ is a product of $e_j$'s only.
\end{Lm}

\begin{proof}  This follows from  Definition \ref{def:aB}
and the fact that there are no nodes lowering $B$.  Indeed,
for $B$ as in the hypotheses, Case (ii) never applies as $\het(B) = 0$,
and it is
immediate in Cases (i) and (iii).
\end{proof}

\begin{Prop}\label{prop:irreducible}
For each simply laced spherical Coxeter type $M$, the algebra $\TL(M)$ has
the following properties.
\begin{enumerate}[(i)]
\item It is  isomorphic to the Temperley--Lieb algebra of type
$M$.
\item The submonoid of $\BrM(M)$ of all monomials
in
$\TL(M)$ (i.e., of height zero) leaves invariant
the collection of all admissible sets in $\AO$ of height zero.
\item For each $Y\in \MY$, the algebra $\TL(M)$ has an
irreducible representation of degree $ |(WB_Y)^0|$.
\item
Up to powers of $\delta$, each monomial
$x$ of $\TL(M)$ is uniquely determined by $ x\emptyset$ and $\emptyset
x$.
\end{enumerate}
\end{Prop}

\begin{proof}
These results are known for $M=\ddA_n$ and $M=\ddD_n$ and only need
to be considered for irreducible Coxeter types, so we restrict attention to
$M=\ddE_6$, $\ddE_7$, $\ddE_8$.

\nl(i). By Lemma \ref{height0aB}, the set $(WB_Y)^0$ is contained in the
orbit of $B_Y$ under $\TL(M)$ in $\AO$. Counting the elements in a monomial
basis of $\TL(M)$ by use of \cite[Lemma 1.3]{CFW}, we conclude that the
rank of $\TL(M)$ is at least
$$\sum_{Y\in \MY} |(WB_Y)^0|^2,$$
which can be seen from Table \ref{table:types2} to be
\begin{eqnarray*}
1+6^2+20^2+15^2 &=& 662,\\
1+7^2+27^2+21^2+35^2+15^2 &=& 2670,\\
1+8^2+35^2+84^2+50^2 &=& 10846,
\end{eqnarray*}
in the respective cases $M = \ddE_6$, $\ddE_7$, $\ddE_8$. These numbers
coincide with the ranks of Temperley-Lieb algebra of type $M$ as
computed by K.~Fan \cite[Section 6.4]{Fan}.
As $\TL(M)$ is a quotient of the Temperley--Lieb algebra of type $M$, we
conclude that it is isomorphic to this Temperley--Lieb algebra.

\nl(ii).  By the equality in (i), the action of each $e_i$ on an element
$B\in (WB_Y)^0$ should stay within $(WB_Y)^0$, for otherwise there would be
too many images of $\emptyset$ in $\AO$ under the monomials in $\TL(M)$ with
regard to (i).

\nl(iii).
Let $Y\in \MY$ and put $\B = WB_Y$.
The restriction to $\TL(M)$ of the linear representation
$\rho_{\B}\otimes1$ of $\Br(M)$ of \cite[Theorem 3.6(ii)]{CFW} is an
irreducible representation of degree $ |(WB_Y)^0|$.
The proof is similar to the proof in \cite[Section~$5$]{CFW}.  Here the
vector space is the linear subspace of $\Br(M)\otimes \Q(\delta)$ with basis
the elements $a_B$ for $B\in (WB_Y)^0$.  To see that this representation is
irreducible, assume $u$ is a nonzero vector in a $\TL(M)$-invariant
subspace.  If $B$ is such that $a_B$ occurs in $u$ with a nonzero
coefficient, we act by $\aback_B$ on $u$ so that the coefficient of $\he_Y =
a_{B_Y}$ is nonzero.  So, without loss of generality, we may assume $\he_Y$
occurs in $u$ with coefficient $1$.  Now multiply $u$ by $e_Y$.  As in
\cite[Proposition~$5.3$]{CFW} all the terms become $\he_Y$ together with a
power of $\delta$. The power of $\delta$ in the coefficient of $\he_Y$ after
this multiplication by $e_Y$ is $\delta^{|Y|}$ and the coefficient of each
other term is a smaller power of $\delta$.  This means $\he_Y$ occurs in the
proper subspace.  But clearly, the span of the images of $\he_Y$ under
$\TL(M)$ is the whole vector space, and so the representation is
irreducible.


\nl(iv).  The map from the basis of Temperley-Lieb monomials
to ordered pairs from $\AO$ in the same $W$-orbit and of
height $0$ given by $x\mapsto (x\emptyset,\emptyset x)$ is well defined by
(ii) and surjective. By (i), $\rk(\TL(M))$ coincides with this number, so the
map is injective as well.
\end{proof}


\begin{Proof} {\bf Proof of Proposition \ref{templieb}}.\label{proof1.2}
\rm By Proposition \ref{prop:irreducible}(i), $\TL(M)$ is the free algebra
with identity generated by $e_i$ $(i=1,\ldots,n)$ subject to the relations
(HSee), (HCee), and (HNeee).  All these relations are homogeneous.  In
particular, the $\isog$-equivalence classes in $F$ having words of height
$0$ correspond bijectively to monomials in $\TL(M)$.  After selecting a
representative for each set of multiples by powers of $\delta$ and extending
the set thus obtained to a set $T$ of reduced words in $F$ such that
$\pi(T)$ is a basis of $\Br(M)$, we can apply Proposition
\ref{prop:BMWbasis}.  This gives us a set $T_0$ of $\rk(\TL(M))$ words in
$\tF$ such that $\rho(T_0)$ is a basis of the subalgebra of $\BMW(M)$
generated by $e_1,\ldots,e_n$. This proves that the subalgebra is isomorphic
to $\TL(M)$, and hence, by Proposition \ref{prop:irreducible}(i) again,
isomorphic to the Temperley--Lieb algebra of type $M$, establishing
Proposition~\ref{templieb}.
\end{Proof}

\begin{Cor} \label{cor:aBszero}
For $B\in\AO$ of height zero,
Theorems \ref{prop:aB} and \ref{th:aBcharacterization} hold.
\end{Cor}

\begin{proof}
Let $Y\in \MY$ and  suppose $B\in (WB_Y)^0$.

We start with Theorem~\ref{prop:aB}.  According to Definition~\ref{def:aB},
the word $a_B$ in $\tF$ has height zero and so its image in $\Br(M)$
belongs to $\TL(M)$.  Also, $\pi(a_B)\emptyset=B$ and
$\emptyset\pi(a_B)=B_Y$ because it is in the same $W$-orbit as $B_Y$ and
$a_B$ ends in $\he_Y$, so $\emptyset\pi(a_B)$ contains the admissible
closure of the set of simple roots indexed by $Y$.  By Proposition
\ref{prop:irreducible}(iv), this determines $a_B$ up to powers of $\delta$.
As existence of $a_B^{\rm b}$ was established in Proposition \ref{property:aB},
Theorem~\ref{prop:aB} follows for $B$ of height zero.

\medskip
We now consider Theorem~\ref{th:aBcharacterization}.  As
$\het(B) = 0$, we have $a_B=e_{i_1}e_{i_2} \cdots e_{i_t}\he_Y$ for certain
nodes $i_1,\ldots,i_t$ of $M$.

Suppose first $|e_iB|>|B|$.  Then $(\a_i,B)=0$, so $\pi(e_ia_B)\emptyset =
\cl{(B\cup\{\a_i\})}$.  Hence there is a set $U\in\MY$ strictly containing
$Y$ such that $ \pi(e_ia_B)\emptyset \in WB_U$.  By Proposition
\ref{prop:irreducible}(ii), the height of $\pi(e_ia_B)\emptyset $ is zero,
and so $\het(e_ia_B)=0 $, from which we conclude that $\emptyset\pi(e_ia_B)
$ has height zero.  But then, by Lemma \ref{prop:irreducible}(iv) applied to
$U$ with admissible sets $\pi(e_ia_B)\emptyset$ and $\emptyset \pi(e_ia_B)$,
respectively, there are elements $a,b\in\TL(M)$ such that $e_ia_B = a \he_U
b^{\op} \in \TL(M)e_U\TL(M)$. This proves (iii).

Suppose then $|e_iB|=|B|$.  Then $\pi(e_ia_B)\emptyset = e_iB\in WB_Y$.
As $\emptyset\pi( e_ia_B)\supseteq B_Y$, we obtain
$\emptyset\pi( e_ia_B)= B_Y$. By Proposition \ref{prop:irreducible}(iv),
this implies  $e_ia_B=a_{e_iB}\delta^k$ for some $k\in\Z$, whence
$e_ia_B\homog a_{e_iB} h$, with $h =\he_Y \delta^k\in H_Y$.
As $\het(e_iB)= \het(B)=0$ has been shown in  Proposition
\ref{prop:irreducible}(ii), we conclude
$\het(e_iB)\le \het(B)$, proving
part (ii) of Theorem~\ref{th:aBcharacterization}.

Finally, we consider Theorem~\ref{th:aBcharacterization}(i).
As $\het (B)=0$, there are no lower elements, so
either $r_iB=B$ or $r_iB>B$. Suppose $r_iB>B$.  Then, by Definition
\ref{def:aB}(ii), $a_{r_iB}=r_ia_B$ and the result follows.

It remains to consider $r_iB=B$.  As $\het(B)=0$, Definition
\ref{def:aB}(iii) applies and gives that $a_B=e_{i_1}e_{i_2}\cdots
e_{i_s}\he_{B_Y}$ for certain nodes $i_1,\ldots,i_s$ of $M$.  We proceed by
induction on the number of terms $e_{i_j}$, which we have denoted $s$, and
prove (i) with $h\in\delta^{\Z}\he_Y$.  If $s=0$, then $B = B_Y$ and $\a_i$
is perpendicular to the simple roots in $B_Y$ or one of these, so
$r_i\he_{Y}\in H_Y$ (for the former case, observe that $H_Y$ contains all
$e_j\he_Y$ with $j\not\sim t$ for all $t\in Y$ and for the latter case, use
(RSre)) and $r_i\he_{Y}\isog\he_{Y}r_i$ (use (HCer) for the former case and
use (RSre) and (RSer) for the latter case).  If $s>0$, Lemma
\ref{sameheight} gives that $B$ contains $\a_{i_1}$ and so $\a_i \perp
\a_{i_1}$ or $\a_i =\a_{i_1}$.  This implies $r_ie_{i_1}e_{i_2}\cdots
e_{i_s}\he_{Y} \isog e_{i_1}r_ie_{i_2}\cdots e_{i_s}\he_{Y}$.  We apply the
induction hypothesis to $B' = e_{i_2}\cdots e_{i_s}B_{Y}$ as $a_{B'}$ has
fewer terms $e_{i_j}$. This gives an exponent $k\in\Z$ such that $r_ia_B\isog
e_{i_1}r_ia_{B'}\homog e_{i_1}a_{r_iB'}\delta^k$.  As
$|e_{i_1}{r_iB'}|=|B|=|r_iB'|$, part (ii) gives $e_{i_1}a_{r_iB'}\homog
a_{e_{i_1}r_iB'}\he_Y\delta^j=a_{r_iB}\he_Y\delta^j$ for some $j\in\Z$, and (i)
follows.
\end{proof}

\end{section}

\begin{section}{Centralizers}\label{sec:WC}
In this section, we establish the rewrite rules for the part of the Brauer
monoid corresponding to the Coxeter group $W(M_{Y})$ as described in
Theorem~\ref{th:nottheta}. This part is the subsemigroup $H_Y$ of Definition
\ref{df:HY}, which centralizes $\he_Y$ in $\tF$.  It
will be shown that the subset $S_Y$ of $H_Y$
is a set of simple reflections of $H_Y$.

Also, we will need $H_{\{n\}}$ to describe a bigger part, to be called
$Z_n$, of the centralizer in $\Br(M)$ of $e_n$. The last result of this
section states that this algebra is a quotient of a Brauer algebra of type
strictly contained in $M$.  These centralizers will help to prove our main
theorems by induction on the rank $n$ of the Coxeter diagram $M$.

\begin{Proof} {\bf Proof of Theorem \ref{th:WC}.} \label{proofWC}
\rm
Let $M = \E_6$ and $Y = \{6\}$. This case corresponds to the
first row of Table \ref{table:types2} below its header. The elements of $S_Y$
are $s_0 = e_6e_5e_4r_2e_3e_4e_5\he_6$ and $s_i = r_i\he_6$ for
$i=1,\ldots,4$.  We have

\begin{eqnarray*}
s_0^2  &= & e_6e_5e_4r_2e_3e_4e_5\he_6e_6e_5e_4r_2e_3e_4e_5\he_6 \homog    e_6e_5e_4r_2e_3e_4e_5e_6e_5e_4r_2e_3e_4e_5\he_6     \\
        &\homog &   e_6e_5e_4r_2e_3e_4e_5e_4r_2e_3e_4e_5e_6 \delta^{-1}  \homog   e_6e_5e_4r_2e_3e_4e_3r_2e_4e_5\he_6    \\
          &\homog &   e_6e_5e_4r_2e_3r_2e_4e_5\he_6     \homog  e_6e_5e_4r_2^2e_3e_4e_5\he_6  \\
          &\homog &    e_6e_5e_4e_3e_4e_5\he_6 \isog  e_6e_5e_4e_5\he_6 \isog \he_6  . \\
\end{eqnarray*}

We next verify the rule $s_1s_0s_1 \homog s_0s_1s_0$.
We are using here that $e_4e_3r_2e_4\isog e_4r_3e_2e_4$ by (HTeere).
\begin{eqnarray*}
s_1s_0s_1&=& r_1e_6e_6e_5e_4e_3r_2e_4e_5e_6r_1e_6\delta^{-3}
 \isog   r_1e_6e_5e_4r_3e_2e_4e_5r_1e_6\delta^{-1} \\
& \isog&  e_6e_5e_4e_2r_1r_3r_1e_4e_5e_6 \delta^{-1}
\isog  e_6e_5e_4e_2r_3r_1r_3e_4e_5\he_6  \\
 &\isog&  e_6e_5e_4r_3e_2r_1r_3e_4e_5\he_6
 \isog  e_6e_5e_4e_3r_4e_2r_1r_4e_3e_4e_5\he_6  \\
 &\isog&  e_6e_5e_4e_3r_4e_2r_4r_1e_3e_4e_5\he_6
\isog  e_6e_5e_4e_3r_2e_4r_2r_1e_3e_4e_5\he_6  \\
 &\isog&  e_6e_5e_4e_3r_2e_4r_1e_3r_2e_4e_5\he_6
\isog  e_6e_5e_4r_2e_3e_4r_1e_3r_2e_4e_5\he_6  \\
 &\isog&  e_6e_5e_4r_2e_3e_4e_5e_6e_5e_4r_1e_3r_2e_4e_5\he_6  \\
&\isog & e_6e_5e_4r_2e_3e_4e_5e_6r_1\he_6e_5e_4e_3r_2e_4e_5\he_6  \\
 &\isog&  e_6e_5e_4r_2e_3e_4e_5\he_6 r_1\he_6e_6e_5e_4e_3r_2e_4e_5\he_6
=  s_0s_1s_0.
\end{eqnarray*}

We next verify the rule $s_2s_0 \isog s_0s_2$.
\begin{eqnarray*}
s_2s_0&=& r_2\he_6e_6e_5e_4e_3r_2e_4e_5\he_6
 \isog   r_2e_6e_5e_4r_3e_2e_4e_5e_6\delta^{-2} \\
& \isog & e_6e_5 r_2e_4e_2r_3e_4e_5e_6\delta^{-2}
 \isog  e_6e_5 r_4e_2r_3e_4e_5e_6\delta^{-2} \\
& \isog&  e_6e_5 e_4r_5e_2r_3e_4e_5e_6\delta^{-2}
 \isog  e_6e_5 e_4e_2r_3r_5e_4e_5e_6\delta^{-2} \\
& \isog&  e_6e_5 e_4r_3e_2r_4e_5e_6\delta^{-2}
 \isog  e_6e_5 e_4e_3r_4e_2r_4e_5e_6\delta^{-2} \\
& \isog&  e_6e_5 e_4e_3r_2e_4r_2e_5e_6\delta^{-2}
 \isog  e_6e_5 e_4e_3r_2e_4e_5\he_6r_2\he_6 \\
&=&  s_0s_2.
\end{eqnarray*}

This settles the case $(M,Y) = (\E_6,\{6\})$.

If $M=\E_6$ and $Y = \{4,6\}$, the elements of $S_Y$ are
$t_1 = r_1\he_{4,6}$ and $t_2 = e_4r_2e_3\he_{4,6}$.
Then obviously $t_1^2\homog \he_{4,6}$ and
\begin{eqnarray*}
t_2^2&=&  e_4r_2e_3e_{4,6}e_4r_2e_3e_{4,6}\delta^{-4}
 \homog e_4e_3r_2e_4r_2e_3e_{4,6}\delta^{-2}\\
 &\homog& e_4e_3r_4e_2r_4e_3\he_{4,6}
\homog e_4r_3e_2r_3\he_{4,6}\\
 &\homog& e_4e_2r_3^2\he_{4,6}
\homog e_4e_2\he_{4,6} \homog \he_{4,6}.
\end{eqnarray*}
Moreover,
\begin{eqnarray*}
t_1t_2t_1&\isog&  r_1e_4e_6e_4r_2e_3e_4e_6r_1e_4e_6\delta^{-6}
\isog  r_1e_4e_2r_3e_4r_1e_4e_6\delta^{-3}  \\
 &\isog& e_4e_2 r_1r_3r_1e_4e_6\delta^{-2}
\isog e_4e_2 r_3r_1r_3e_4e_6\delta^{-2}  \\
 &\isog& e_4 r_3r_1e_2r_3e_4e_6\delta^{-2}
\isog e_4 e_3r_4 r_1e_2r_4e_3e_4e_6  \delta^{-2}\\
 &\isog& e_4 e_3 r_1 r_4e_2r_4e_3e_4e_6  \delta^{-2}
\isog e_4 e_3 r_1 r_2e_4r_2e_3e_4e_6  \delta^{-2}\\
 &\isog&  e_4r_2e_3e_4r_1e_4r_2e_3e_4e_6 \delta^{-3}
\isog  e_4r_2e_3e_4e_6r_1e_4e_6e_4r_2e_3e_4e_6 \delta^{-6} \\
 &\isog&  t_2t_1t_2.  \\
\end{eqnarray*}
This settles the case $(M,Y) = (\E_6,\{4,6\})$.  For the case $(M,Y) =
(\E_6,\{2,3,5\})$ there is nothing to prove except $\he_2\he_3\he_5$ is an
idempotent, which follows as $\{2,3,5\}$ is a coclique in $M$.  This settles
the first part of Theorem \ref{th:WC} on the Matsumoto--Tits rules. The
second part on the bijective correspondence follows as the image $\pi(H_Y)$
is known to be of size $W(M_Y)$ from \cite[Lemma 1.3]{CFW}.

Similar computations work for $M = \E_7$ and $M = \E_8$.
\end{Proof}

We derive the following consequence, in which $l$ is the usual length
function on Coxeter groups.

\begin{Cor}\label{Cor:WC}
Let $Y\in\MY$.  The map $M_Y \to H_Y$ sending the $i$-th simple reflection
of the Weyl group $W(M_Y)$ to the $i$-th element listed in the column of
Table \ref{table:types2} for $S_Y$ induces an isomorphism of Coxeter groups
$\zeta_Y: W(M_Y)\to \pi(H_Y) $. In particular, for $w\in W(M_Y)$, we have
$l(w) = \het(\zeta_Y(w))$.
\end{Cor}

\begin{proof}  Theorem~\ref{th:WC} gives a surjective homomorphism of
monoids.  We use \cite[Proposition~$4.7$ (iii)]{CFW}.  We use here $B_Y$
of the table rather than the highest element of the poset $WB_Y$ as in
\cite{CFW}.  By \cite[Lemma~$4.4$]{CFW} we see $e_X$ there corresponds to
$e_Y$ here up to a power of $\delta$.  We have chosen the elements of $S_Y$
to be the generators of the complement to $A_X$ multiplied by $\he_Y$ in
\cite[Proposition~$4.7$ (iii)]{CFW}.  This means the size of $\pi(H_Y)$
coincides with $|W(M_Y)|$, so it is an isomorphism of monoids. As $W(M_Y)$ is a
group, it is an isomorphism of groups as well.

Note that the generators we have chosen in $S_Y$ all have height one as do
the generators of $W(M_Y)$ and so $l(w)=\het( \zeta_Y(w))$ for each $w\in
W(M_Y)$.
\end{proof}

\np Theorem~\ref{th:WC} exhibits a subsemigroup of $\BrM(M)$ isomorphic to
the Coxeter group $W(M_Y)$ for the particular case $Y = \{n\}$.  We
introduce the word $f_0=e_ne_{n-1}\cdots e_4e_2e_3e_4 \cdots e_{n-1}\he_n$
and, for each $i$ with $1\leq i \leq n-2$, the word $f_i=e_i\he_n$ in $F$.
In other words, the $f_i$ are the same as the $s_i$ for $Y = \{n\}$ of Table
\ref{table:types2}, but with the single $r_2$ that occurs in their defining
expression replaced by $e_2$. Now $Z_n$ is defined as the nonunital
subalgebra of $\Br(M)$ generated by $\pi(H_{\{n\}})$ and the images of
$f_0,f_1,\ldots ,f_{n-2}$ under $\pi$; then $Z_n$ has identity element
$\he_n$.  We will extend the group homomorphism $\zeta_{\{n\}}:
W(M_{\{n\}})\to \pi(H_{\{n\}})$ of Corollary~\ref{Cor:WC} to a surjective
algebra homomorphism $\Br(M_{\{n\}}\to Z_n$ for the cases $M = \E_n$ where
$n=6,7,8$.  (Recall that $\Br(M_Y)$ is the algebra generated by the
generators and relations of Table~\ref{BrauerTable}.)  For ease of
presentation, we will write $H_n$, $M_n$, and $\zeta_n$, instead of
$H_{\{n\}}$, $M_{\{n\}}$, and $\zeta_{\{n\}}$, respectively.  Clearly, the
subalgebra $Z_n$ contains $\pi(H_n)$ and has identity element $\he_n$.

\begin{Prop}\label{th:brauersubalgebra}
Let $n\in \{6,7,8\}$ and $M=\E_n$. Take $Y = \{n\}$ and consider the
diagram $M_Y = M_n = \A_5$, $\D_6$, $\E_7$ if $n=6,7,8$,
respectively. The rewrite rules of Table
\ref{BrauerTable}
with respect to $\homog$ for type $M_n$ are satisfied by
$s_0,s_1,\ldots,s_{n-2}$ instead of the $r_i$ and $f_0,f_1,\ldots,f_{n-2}$
instead of the $e_i$.
In particular, there is a surjective algebra homomorphism $\zeta_n:
\Br(M_n)\to Z_n$ determined by $\zeta_n(r_i) = s_i$ and $\zeta_n(e_i)
= f_i$, for $0\le i\le n-2$.
\end{Prop}

Here the labeling for $M_n$ is as in the subdiagram of
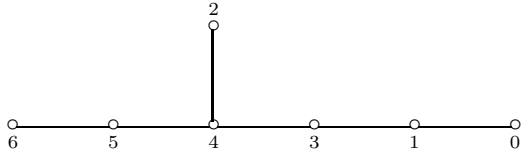
\begin{eqnarray*}
&&\Esevensigma
\end{eqnarray*}
induced on $\{0,\ldots,n-2\}$.  So the full diagram is for $\E_8$; for
$\E_7$, delete $6$; for $\E_6$, delete $6$ and $5$.

\begin{proof}
We treat the case $n=6$ and leave the other cases to the reader.
We check that the powers of $\delta$ work as required.  In view of
Theorem~\ref{th:WC},
the only new relations needed are the ones involving $f_i$.  These are
all straightforward unless one of the indices is $0$.
For instance,
if $i\neq 0$, then
$f_i^2=e_i^2 e_6^2\delta^{-2}\isog e_i e_6 =\delta e_i\he_6 = \delta f_i$.
Moreover,
\begin{eqnarray*}
f_0^2  &\homog & e_6e_5e_4e_2e_3e_4e_5e_6e_6e_5e_4e_2e_3e_4e_5\he_6\delta^{-1} \homog    e_6e_5e_4e_2e_3e_4e_5e_6e_5e_4e_2e_3e_4e_5\he_6     \\
        &\homog &   e_6e_5e_4e_2e_3e_4e_5e_4e_2e_3e_4e_5\he_6  \homog   e_6e_5e_4e_2e_3e_4e_3e_2e_4e_5\he_6    \\
          &\homog &   e_6e_5e_4e_2e_3e_2e_4e_5\he_6     \homog  e_6e_5e_4e_2^2e_3e_4e_5\he_6  \\
          &\homog &    e_6e_5e_4e_2e_3e_4e_5\he_6\delta \isog  \delta f_0 ,
\end{eqnarray*}
and so (HSee) is satisfied.

These same equations are easily modified to verify (RSre) and (RSer) for the
cases $s_0$ and $f_0$.  In particular we need $s_0f_0=f_0s_0=f_0$.  As for
$s_0f_0$, the leftmost $e_2$ in the above reduction for $f_0^2$ becomes
$r_2$ in the definition of $s_0$.  Follow the equations using the same
relations until the occurrence of $e_2^2$, which becomes $r_2e_2$ and so
reduces to $e_2$.  The result follows (without the appearance of $\delta$).

We verify the instance $s_1f_0s_1 \isog s_0f_1s_0$ of (HNrer).
\begin{eqnarray*}
s_1f_0s_1&=& r_1\he_6e_6e_5e_4e_3e_2e_4e_5\he_6r_1\he_6
 \isog   r_1e_6e_5e_4e_3e_2e_4e_5r_1\he_6 \\
& \isog&  e_6e_5e_4e_2r_1e_3r_1e_4e_5\he_6
\isog  e_6e_5e_4e_2r_3e_1r_3e_4e_5\he_6  \\
 &\isog&  e_6e_5e_4r_3e_2e_1r_3e_4e_5\he_6
 \isog  e_6e_5e_4e_3r_4e_2e_1r_4e_3e_4e_5\he_6  \\
 &\isog&  e_6e_5e_4e_3r_4e_2r_4e_1e_3e_4e_5\he_6
\isog  e_6e_5e_4e_3r_2e_4r_2e_1e_3e_4e_5\he_6  \\
 &\isog&  e_6e_5e_4e_3r_2e_4e_1e_3r_2e_4e_5\he_6
\isog  e_6e_5e_4r_2e_3e_4e_1e_3r_2e_4e_5\he_6  \\
 &\isog&  e_6e_5e_4r_2e_3e_4e_5e_6e_5e_4e_1e_3r_2e_4e_5\he_6  \\
&\isog & e_6e_5e_4r_2e_3e_4e_5e_6e_1\he_6e_5e_4e_3r_2e_4e_5\he_6  \\
 &\isog&  e_6e_5e_4r_2e_3e_4e_5\he_6 e_1\he_6e_6e_5e_4e_3r_2e_4e_5\he_6
=  s_0e_1s_0.
\end{eqnarray*}

We next derive the instance $s_2f_0 \isog f_0s_2$ of (HCer).
\begin{eqnarray*}
s_2f_0&=& r_2\he_6e_6e_5e_4e_3e_2e_4e_5\he_6
 \isog   r_2e_6e_5e_4e_3e_2e_4e_5e_6\delta^{-1} \\
& \isog & e_6e_5 r_2e_4e_2e_3e_4e_5e_6\delta^{-1}
 \isog  e_6e_5 r_4e_2e_3e_4e_5e_6\delta^{-1} \\
& \isog&  e_6e_5 e_4r_5e_2e_3e_4e_5e_6\delta^{-1}
 \isog  e_6e_5 e_4e_2e_3r_5e_4e_5e_6\delta^{-1} \\
& \isog&  e_6e_5 e_4e_3e_2r_4e_5e_6\delta^{-1}
 \isog  e_6e_5 e_4e_3e_2e_4r_2e_5e_6\delta^{-1} \\
& \isog&  e_6e_5 e_4e_3e_2e_4e_5e_6r_2\delta^{-1}
  \isog  e_6e_5 e_4e_3e_2e_4e_5e_6\he_6r_2\he_6
=  f_0s_2.
\end{eqnarray*}

Now we consider $e_2f_0 \isog f_0e_2$; we have
\begin{eqnarray*}
e_2f_0&=& e_2\he_6e_6e_5e_4e_3e_2e_4e_5\he_6
 \isog   e_2e_6e_5e_4e_3e_2e_4e_5e_6\delta^{-1} \\
& \isog & e_6e_5 e_2e_4e_2e_3e_4e_5e_6\delta^{-1}
 \isog  e_6e_5 e_2e_3e_4e_5e_6\delta^{-1}. \\
\end{eqnarray*}
This is symmetric (fixed under $\op$) as $e_2$ and $e_3$ commute and so is
homogeneously equivalent to $f_0e_2$.

The remaining rewrite rules are easily verified in the same manner. We only
treat (RNrre) here.  There are two instances involving $s_0$.  First there
is $s_1s_0f_1\homog f_0 f_1$, which we verify as follows.
\begin{eqnarray*}
s_1s_0f_1 &=& r_1\he_6 e_6e_5e_4r_2e_3e_4e_5\he_6e_1 \he_6 \isog r_1e_6e_5e_4r_2e_3e_4e_5e_6e_1\delta^{-1} \\
          &\isog & e_6e_5e_4r_2r_1e_3e_1e_4e_5e_6\delta^{-1} \homog e_6e_5e_4r_2r_3e_1e_4e_5e_6\delta^{-1}  \\
          &\homog& e_6e_5e_4e_2e_3e_4e_5e_6e_1\delta^{-1}     =      f_0 f_1.
\end{eqnarray*}

To finish this, we need to verify is $s_0s_1f_0\isog f_1f_0$.
\begin{eqnarray*}
s_0s_1f_0  &=& e_6e_5e_4r_2e_3e_4e_5\he_6 r_1 \he_6 e_6e_5e_4e_2e_3e_4e_5\he_6      \\
                &\homog&       e_6e_5e_4r_2e_3e_4e_5e_6 r_1  e_6e_6e_5e_4e_2e_3e_4e_5e_6 \delta^{3}  \\
                        &\homog&  e_6e_5e_4r_2e_3e_4e_5e_6 r_1 e_5e_4e_2e_3e_4e_5e_6 \delta  \\
                 &\homog&  e_6e_5e_4e_3r_2r_1e_4e_5e_6  e_5e_4e_3e_2e_4e_5e_6 \delta    \\
                 &\homog &   e_6e_5e_4e_3r_2r_1e_4e_5e_6  e_5e_4e_3e_2e_4e_5e_6 \delta
\homog            e_6e_5e_4e_3r_2r_1e_4 e_3e_2e_4e_5e_6 \delta    \\
 &\homog &            e_6e_5e_4e_3r_2e_4e_3r_1 e_3e_4e_5e_6 \delta
\homog            e_6e_5e_4e_3r_2e_4e_2r_1 e_3e_4e_5e_6 \delta    \\
 &\homog &            e_6e_5e_4e_3r_2e_2r_1 e_3e_4e_5e_6 \delta
\homog   e_6e_5e_4e_3r_4e_2r_1 e_3e_4e_5e_6 \delta    \\
 &\homog &            e_6e_5e_4r_3e_2r_1 e_3e_4e_5e_6 \delta
\homog             e_6e_5e_4e_2r_3r_1 e_3e_4e_5e_6 \delta    \\
 &\homog &            e_6e_5e_4e_2e_1e_3e_4e_5e_6 \delta
\homog             e_1e_6e_5e_4e_2e_3e_4e_5e_6 \delta    \\
   &\homog &      e_1\he_6 e_6e_5e_4e_3e_2e_4e_5\he_6 \homog f_1f_0.
\end{eqnarray*}
\end{proof}

\begin{Remark}\label{th:notfullbrauersubalgebra}
\rm According to Proposition~\ref{th:brauersubalgebra}, the algebra $Z_n$ is
a homomorphic image of $\Br(M_n)$.  Unlike many of the properties of
subalgebras generated by subsets of the generators, $Z_n$ is not the full
Brauer algebra, but is a proper quotient.
We will show this for
$n=6$ by exhibiting two distinct
elements in $\Br(M_6)$ whose images are the same in $Z_6$.
Recall that $\Br(M_{6})$ has type $\A_5$.  The fundamental roots of $M_{6}$
can be taken to be $\{\a_2,\a_4,\a_3,\a_1,\a_0\}$ with
$\a_0=\a_2+\a_3+2\a_4+2\a_5+\a_6$ which is the highest root of the root
system of type $\D_5$ spanned by $\a_i$ for $i\geq 2$ within the root system
of type $\E_6$.
The elements $e_2e_3$ and $e_2e_3e_0$ are distinct in $\Br(\A_5)$ (with the
labeling as in the above diagram for $M_n$), but their $\zeta_6$-images
$\pi(f_2f_3)$ and $\pi(f_2f_3f_0)$ coincide in $Z_6$, as
$e_2e_3(e_6e_5e_4e_3e_2e_4e_5e_6)\isog\delta e_2e_3e_6$ (obtained by
straightforward reductions).  These elements are not $0$ in $\Br(\E_6)$
by the results of \cite{CFW}.
Therefore, $Z_6$ is a proper quotient of $\Br(\A_5)$.

The same ideas work for $n=7$ and $8$.
\end{Remark}


\np
The image of $\BrM(M_n)$ in $Z_n$ under $\zeta_n$
of Proposition \ref{th:brauersubalgebra} is a monoid acting on
$\AO$, and so we can view the monoid $\BrM(M_n)$ itself as acting on
$\AO$.  For a subset $\B$ of $\AO$, denote by $\B^n$ the set of those
admissible sets in $\B$ that contain $\a_n$, and
by $\B^*$ the set of all $B'\setminus\{ \a_n\}$ for $B'\in\B^n$.

\begin{Lm}\label{Bsetminusan}
The set $\AO^*$ consists of admissible sets for $\BrM(M_n)$.
If $\B$ is a $W(M)$-orbit in $\AO$, then
$\B^*$ is a $W(M_n)$-orbit in $\AO^*$.
\end{Lm}

\begin{proof}
Let $B^*\in \AO^*$, so $B = B^*\cup\{\a_n\}\in\AO$.
As the elements of $B$ are mutually orthogonal, so are the elements of
$B^*$.  The action of a reflection from $\zeta_n(W(M_n))$ on
the set $B$ fixes $\a_n$ and, because $B$ is admissible, the reflection
moves $0,1,2,4$ points by \cite[Proposition~$2.3$, (iii)]{CGW2};
consequently it moves the same number of points in $B^* = B\setminus \{\a_n\}$.
Now by this same proposition, $B^*$ is admissible.

The group $W(M_n)$ is a submonoid of $\BrM(M_n)$ and so acts on $\AO$ via
$\zeta_n$. Each of its elements fixes $\a_n$. Therefore, $W(M_n)$ leaves
$\B^n$ invariant, and hence also $\B^*$.  To show $W(M_n)$ is transitive on
$\B^*$, we consider two elements $B'$ and $B''$ of $\B^n$.  As they are in
the same $W$-orbit, there is an element $w\in W$ with $wB'=B''$.  For
each such $B'$ the action of the normalizer in $W$ of $B'$ is given in
\cite[Table~$3$]{CFW} and in each case, it is transitive on $B'$.  We can
then act by an element of the normalizer to ensure that $w$ takes $\a_n$ to
$\a_n$.  This implies $w\in W(M_n)$ by a well-known result on reflection
groups (\cite[Exercice V.6.8]{Bour}).  As $w$ takes $B'\setminus\{\a_n\}$
to $B''\setminus\{\a_n\}$, we conclude that $W(M_n)$ is transitive on
$\B^*$.
\end{proof}

\np A look at Table \ref{table:types2} shows that, for $M$ of type $\E_n$,
the $W$-orbits in $\AO$ are uniquely determined by the size of a
representative element. This is not the case for $M = \D_n$.  For each
$W$-orbit $\B$ of admissible sets of given size $k$, except for $M=\E_7$
with $k=3$ or $4$, there is a unique $W(M_n)$-orbit of admissible sets
of size $k-1$, so $\B^*$ is uniquely determined by $k$.  In the case where
$M = \E_7$, we have $M_n=\D_6$ and there are two $W(\D_6)$-orbits of
admissible sets of size $3$.  Here, the $W(\D_6)$-orbit arising as $\B^*$
from $\B$ for $k=4$ is the one containing $\{\a_0, \a_3, \a_2\}$ where
$\a_0$ is the root $\a_2+\a_3+2\a_4+2\a_5+2\a_6+\a_7$, rather than
$\{\a_0, \a_3, \a_5\}$.  This can be seen by starting with $B =
\{\a_3,\a_2,\a_5,\a_3+\a_2+2\a_4+\a_5\}$ and acting by $r_6r_7r_5r_6$.  For
$k=3$, the admissible sets of size two contain $\a_5$ and $\a_2$ and so
$\B^*$ is the orbit of size $15$ in the second line of \cite[Table~3]{CFW}
for $\D_6$.  The sizes are listed in Table~\ref{table:blockswithn} which can
be obtained either directly as indicated here or by using GAP, \cite{GAP}.

\begin{table}[ht]
\begin{center}
\begin{tabular}{|cc|cc|}
\hline
${\CoxDiag}$ & $|B_Y|$ & $M_n$ &$|\B^n|$\\
\hline
$\E_6$ & $2$ & $\A_5$ & $15$ \\
$\E_6$ & $4$ & $\A_5$ & $15$ \\
\hline
$\E_7$ & $2$ & $\D_6$ & $30$ \\
$\E_7$ & $3$ & $\D_6$ &$15$ \\
$\E_7$ & $4$ &  $\D_6$ &$60$ \\
$\E_7$ & $7$ & $\D_6$ & $15$  \\
\hline
$\E_8$ & $2$ & $\E_7$ & $63$ \\
$\E_8$ & $4$ & $\E_7$ & $315$ \\
$\E_8$ & $8$ & $\E_7$ & $135$ \\
\hline
\end{tabular}
\end{center}
\smallskip
\caption{\label{table:blockswithn}\textrm{Numbers $|\B^n|$
of sets in $\B = WB_Y$ containing $\a_n$
for $Y\in\MY$ with $|B_Y|\ge2$.}}
\end{table}

In Corollary~\ref{heightsthesame} we will show that the height of $B$ in
the poset $\AO$ for $M = \E_n$ is the same as the height of $B\setminus\{
\a_n\}$ in the poset for $M_n$ with sets of this size.

\begin{Notation}\label{df:Zj}
\rm
Let $k_1, k_2, \ldots, k_l$ be a sequence of nodes of $M$.  Then
$e_{k_1}e_{k_2}\cdots e_{k_l}$ will be denoted by $e_{k_1,\ldots, k_l}$.  In
the special case where $k=k_1,\ldots,k_l=j$ is the path from $k$ to $j$ in
$M$, we also write $e_{kj}$. Moreover, we adopt the same notation for the
hatted versions, e.g., $\he_{k_1,\ldots,k_l} = \he_{k_1}\cdots \he_{k_l}$.
For $j\in \{1,\ldots,n\}$ we write $Z_j = \he_{jn}Z_n\he_{nj}$.
\end{Notation}

\begin{Lm}\label{lm:Zj}
The algebra $Z_j$ is isomorphic to $Z_n$ via the height preserving maps
$x\mapsto \he_{nj}x\he_{jn}:Z_j\to Z_n$ and $y\mapsto
\he_{jn}y\he_{nj}:Z_n\to Z_j$.  Moreover these algebras satisfy the same
rewrite rules for type $M_n$ as stipulated in Proposition
\ref{th:brauersubalgebra} with respect to their natural generators.
Accordingly, $\AO^j = e_{jn}\AO^n =
\{e_{jn}B\mid B\in\AO^n\}$ is the set of all
admissible elements containing $\a_j$ and satisfies $e_{nj} \AO^j = \AO^n$.
\end{Lm}

\begin{proof}
By (HNeee) $\he_{jn}\he_{nj} = \he_j$ and so the map $y\mapsto
\he_{jn}y\he_{nj}$ on $Z_n$ is the inverse of $x\mapsto \he_{nj}x\he_{jn}$
on $Z_j$.  As $\he_{nj}$ has height $0$ and all $x\in Z_j$ commute with
$\he_j$, the assertions about rewrites follow.

Finally, if $j=j_1,j_2,\ldots,j_l=k$ is the path in $M$ from
$j$ to $k$, then, for $B\in \AO^j$, we have $\a_j\in B$ by Lemma
\ref{sameheight}, so
$\he_{nj}B$ is obtained from $B$ by applying the Howlett-Brink word
$r_{j_{l-1}}r_{j_l} \cdots r_{j_2}r_{j_3}r_{j_1}r_{j_2}\he_j=
e_{j_l, \ldots, j_1} = e_{k j}$. We conclude that $e_{kj}\AO^j = \AO^k$.
\end{proof}

There is an important property that lowering-e-nodes possess.
\begin{Lm}\label{lm:uniquelowenode}
Suppose that $l$ is a lowering-e-node for $B$ and $j\sim l$ satisfies
$\a_j\in B$.  Then $\a_l$ is orthogonal to every simple root in
$B\setminus\{\a_j\}$.
\end{Lm}

\begin{proof}
If $l\sim k$ with $\a_{k}\in B\setminus\{\a_j\}$, then $r_jr_l$ would map
the pair $\{\a_j,\a_{k}\}$ to $\{\a_{l},\a_l+\a_{k} +\a_j\}$, and so the
level of $e_lB=r_jr_lB$ would be higher than $L(B)$, contradicting
$L(e_lB)<L(B)$.
\end{proof}

\begin{Notation}
\label{not:NlB}
\rm By Lemma \ref{lm:uniquelowenode}, for each lowering-e-node $l$ for $B$,
there is a unique simple root $\a_j$ in $B$ such that $j\sim l$, and we
write $j = N(l,B)$.
\end{Notation}

\np
The following lemma exhibits elements of $Z_j$ which appear in the lowering
algorithm of Definition \ref{def:aB}.

\begin{Lm}\label{inZi}
Suppose $\a_j\in B$ and $i_1,\ldots,i_t$ is a string of nodes of $M$ such
that each ${i_j}$ is a lowering-e-node for
$e_{i_{j-1}}\cdots e_{i_2}e_{i_1}B$.
Now set $B^{{\rm i}}=e_{i_t}e_{i_{t-1}} \cdots e_{i_1}B$
and assume $B^{{\rm i}}$ is the first one with a lowering node $s$, so
$B^{{\rm ii}}=r_sB^{{\rm i}} <B^{{\rm i}}$.
For each $k\in\{1,\ldots,t\}$, put $j_k = N(i_{k},e_{i_{k-1}}\cdots e_{i_1}B)$.
Then $e_je_{j_1\cdots j_t}r_se_{i_t\cdots i_1}e_j\in Z_j$.
\end{Lm}

\begin{proof}
We proceed by induction on $t$.

Assume $t=1$. Set $i=i_1$ and $k = j_1$.
If $k=j$, the word under consideration is $e_je_kr_se_ie_j = e_jr_se_ie_j$
(observe that $s\not\sim i$ as $\a_i\in B^{\rm i}$ and $s$ lowers $B^{\rm
i}$), which is in $Z_j$ as $r_se_i\in Z_i$ and $Z_j=e_jZ_ie_j$.

If $k\neq j$ we get $e_je_kr_se_ie_j$.  But
by Lemma \ref{lm:uniquelowenode}, there is only one root in $B$,
namely $\a_k$, not orthogonal to $\a_i$,
so $k\not \sim i$ and $k \not \sim j$ as $\a_k$ and $\a_j $ are in
$B$ and so are orthogonal.  Now $e_ke_i$ and $e_ke_j$ are in $Z_k$.  Also
$s\not \sim k$ (for otherwise $\a_k$ would be raised by $r_s$) and so
$r_se_k\in Z_k$ also.  In particular $e_je_kr_se_ie_j =
e_ke_jr_se_ie_j\in Z_k$.

Suppose then $t>1$.
Now use induction and consider $w=e_{j_2\cdots j_t}r_se_{i_t\cdots i_1}$.
If $k=j=j_1$, then, as $\a_{i_1}\in e_{i_1}B$, by induction
$e_{i_1}w=\delta w\in Z_{i_1}$ and then $e_jwe_j\in Z_j$ as $i_1\sim j=j_1$.
If $k\neq j$ then $k \not \sim i_1$, giving $\a_k\in e_{i_1}B$ and so
by induction $e_kw\in Z_k$.  But then $e_ke_j$ and $e_ke_{i_1}$ are in $Z_k$
finishing the lemma.
\end{proof}

There is an immediate corollary.  Recall the terminology of
Lemma~\ref{Bsetminusan} in which $\AO_n$ is the subset of $\AO$ for which
each set contains $\a_n$ as one of its orthogonal roots and $\AO^*$ is the
set of all $B\setminus \{\a_n\}$ for $B\in \AO^n$.

\begin{Cor}\label{heightsthesame}
Suppose $B\in \AO^n$.  Then the height of $B$ in the poset $\AO$ is the same
as the height of $B^*$ in the poset $\AO^*$.
\end{Cor}

\begin{proof}
Let $Y\in \MY$ and $B\in WB_Y$.  The height of $B$ in the poset for $\E_n$
is the number of terms $r_i$ in $a_B$ by Proposition~\ref{property:aB}.  The
height of $B\setminus\{ \a_n\}$ in the poset $\AO^*$ is the number of
reducing steps it takes to reduce $B\setminus \{\a_n\}$ to a set with
$|Y|-1$ simple nodes.  We know this can be done in $\het(B)$ steps by the
construction above.  These are all lowering moves and so $\het(B)$ is the
height of $B\setminus\{ \a_n\}$ in the poset $\AO^*$.
\end{proof}

\end{section}

\begin{section}{Properties of $a_B$}\label{sec:aB}
This section is devoted to the proof of Theorem \ref{prop:aB}. We fix
$Y\in\MY$ and $B\in WB_Y$.  Throughout the section, we assume the truth of
this theorem and Theorem \ref{th:aBcharacterization} for admissible sets of
level smaller than $B$.

The height zero cases of both theorems were proved in
Corollary~\ref{cor:aBszero}.  Therefore, we can and will assume $\het(B)>0$.
We will also use induction on the rank $n$ of $M$. Recall the validity
of both theorems for simply laced Coxeter diagrams $M$ of type $\A_m$
$(m\ge1)$ and $\D_m$ $(m\ge4)$.

Existence of $a_B$ and $a_B^{\back}$ in $F$ is proved in
Proposition~\ref{property:aB}(iii). For the uniqueness proof, we only need
consider $a_B$; we distinguish the three cases of Definition \ref{def:aB}.

\nl Case (i).  If $|\Simp(B)| = |\Simp(B_Y)|$, then $B$, being the admissible
closure of a set of simple roots, has height~$0$ and so the
statement follows from Corollary \ref{cor:aBszero}.

\nl Case (iii).  Here $|\Simp(B)| < |\Simp(B_Y)|$ and $r_jB\ge B$ for each
node $j$ of $M$.
Then there is a simple root $\a_j$ in $B$.

We will rewrite $a_B$ homogeneously to a product of a monomial in $Z_j$ (see
Definition \ref{df:Zj}) of height $\het(B)$ and a monomial of $\TL(M)$ (see
Notation \ref{not:TL}).

By Definition~\ref{def:aB}
there is a string of nodes $\{i_1,i_2,\ldots , i_t\}$ which are successive
lowering-e-nodes for $B$, $e_{i_1}B$, $e_{i_2}e_{i_1}B$, etc.  Now set
$B^{{\rm i}}=e_{i_t}e_{i_{t-1}} \cdots e_{i_1}B$ and assume $B^{{\rm i}}$ is
the first one with a lowering node $s$. Thus, $B^{{\rm ii}}=r_sB^{{\rm
i}}<B^{{\rm i}}$.
For each $k\in\{1,\ldots,t\}$, put $j_k = N(i_{k},e_{i_{k-1}}\cdots
e_{i_1}B)$
(see Notation \ref{not:NlB}).
By Lemma \ref{inZi},
the monomial $e_je_{j_1\cdots j_t}r_se_{i_t\cdots i_1}e_j$ belongs to
$Z_j$.

By definition, $a_B=e_{j_1\cdots j_t}r_sa_{B^{\rm ii}}$ where $B^{\rm
i}=e_{i_t\cdots i_1}B$ and $B^{\rm ii}=r_sB^{\rm i}$.  Set $B^{{\rm
iii}}=e_{j_1\cdots j_t}B^{{\rm ii}}$.  Then $B^{{\rm iii}}$ also contains
$\a_j$ and so $a_{B^{{\rm iii}}}\isog e_ja_{B^{{\rm iii}}}$.
By induction, $a_{B^{{\rm ii}}} \isog e_{i_t\cdots i_1}a_{B^{{\rm iii}}}$
and so $a_B\isog d a_{B^{{\rm iii}}}$ where $d = e_{j_1\cdots
j_t}r_se_{i_t\cdots i_1}e_j$.  If $j_1=j$, then $d =
\delta^{-1}e_je_{j_1\cdots j_t}r_se_{i_t\cdots i_1}e_j\in Z_j$.
If $j_1\ne j$, then $j\not\sim i_1$ by Lemma \ref{lm:uniquelowenode},
and, by induction and Theorem \ref{th:aBcharacterization}(ii),
as $L(e_{i_1}B)<L(B)$, we have $e_j a_{e_{i_1}B} \homog \delta^p
a_{e_{i_1}B}$ for some integer $p$,
so $a_B = e_{j_1}a_{e_{i_1}B} \isog \delta^{-p} e_{j_1}e_j a_{e_{i_1}B}
\isog \delta^{-p} e_j e_{j_1} a_{e_{i_1}B} \isog \delta^{-p} e_j a_B$.  We
conclude $a_B\isog \delta^{-p}e_ja_B\isog \delta^{-p} e_j d a_{B^{{\rm
iii}}}$ with $e_jd\in Z_j$, so $a_B$ is homogeneously equivalent to
$za_{B^{{\rm iii}}}$, where $z$ is a monomial in $ Z_j$ of height $1$ and
$B^{{\rm iii}}$ contains $\a_j$
and has height
$\het(B) - 1$.

Now compute $a_{B^{{\rm iii}}}$ working only in $Z_j$ and using the set
$\AO^j$ of elements containing $\a_j$ as one of the roots.  By induction on
$M$ the word $a_{B\setminus\{\a_j\}}$ for $M_n$, denoted
$a'_{B\setminus\{\a_j\}}$, is unique up to powers of $\delta$ and homogeneous
equivalence.  Here, the basic height $0$
admissible element for $Z_j$ in $\AO^j$ is taken to be $C = e_{jn}B_Y$.
By Theorem \ref{th:WC}, $a_B$ is homogeneously
equivalent to $a_Ca'_{B\setminus\{\a_j\}}$.
By Corollary \ref{cor:aBszero}, the word $a_C$ is also unique up to powers
of $\delta$ and homogeneous equivalence.
This establishes Case (iii).

\nl Case (ii). Here we use
\cite[Proposition~$3.1$]{CGW2}, \cite[ Lemma~$3.2$]{CGW2} and
\cite[Lemma~$3.3$]{CGW2} which we record here as lemmas for the convenience
of the reader.  We continue to let $\B$ be an $W$-orbit in $\AO$.

\begin{Lm}\label{Prop:AdmOrd}\cite[Proposition $3.1$]{CGW2}
The
ordering $<$ on $\AO$ has the following properties.
\begin{enumerate}[(i)]
\item For each node $i$ of $M$ and each $B\in \B$, the sets $B$
and $r_i B$ are comparable. Furthermore, if $(\a_i,\b)=\pm1$ for some
$\b \in B$, then $r_iB\ne B$.
\item Suppose $i\sim j$ and $\a_i\in
B^\perp$. If $r_jB < B$, then $r_ir_jB<r_jB$. Also, $r_jB > B$
implies $r_ir_jB > r_j B$. \item If $i\not\sim j$, $r_i B<B$,
$r_jB < B$, and $r_i B\ne r_j B$, then $r_ir_jB<r_jB$ and
$r_ir_jB<r_i B$. \item If $i\sim j$, $r_i B<B$, and $r_jB < B$,
then either $r_ir_jB = r_jB$ or $r_ir_jB<r_jB$, $r_jr_iB<r_i B$,
$r_ir_jr_iB< r_ir_jB$, and $r_ir_jr_iB<r_jr_iB$.
\end{enumerate}
\end{Lm}

\begin{Lm}\label{heighttworoot}\cite[Lemma~$3.2$]{CGW2}
Suppose that $B\in\B$ satisfies $r_ir_jB=r_jB$ with $i\sim j$.  If $r_iB<B$
and $r_jB<B$, then $\a_i+\a_j\in B$.
\end{Lm}

\begin{Lm}\label{goingupequal}\cite[Lemma~$3.3$]{CGW2}
Suppose $B\in \B$ and $r_iB=r_kB>B$ with $k\ne i$.  If $\b $ is
the element of $B$ of smallest height moved by either $r_i$ or
$r_k$, then $\b+\a_i+\a_k$ is also in $B$.  Furthermore, $i\not
\sim k$.
\end{Lm}

Assume now that $B$ has two different lowering nodes, $l$ and $k$, so
$r_lB<B$ and $r_kB<B$.  We assume first that $l\not \sim k$.  Using
Lemma~\ref{Prop:AdmOrd}(iii) we see either $r_kB=r_lB $ or both
$r_lr_kB<r_lB$ and $r_lr_kB<r_kB$.  If $r_kB \neq r_lB$, the path down which
starts with $r_l$ can be continued down with $r_k$.  By induction this gives
$a_B \isog r_la_{r_lB}\isog r_lr_ka_{r_kr_lB}$.  Do the same for the path
which starts with $r_k$ and continues with $r_l$; the result
is $a_B \isog r_kr_la_{r_lr_kB}$, which is homogeneously equivalent to the
previous expression.

We next assume $r_lB=r_kB$ (and still $l\not\sim k$).

\begin{Lm}\label{equallowerings}  Suppose $r_lB=r_kB<B$.
Suppose further there is a node $j$ for which $r_jB<B$ with $r_jB\neq B^{\rm
i}$.  Then $r_la_{r_lB}\homog r_ka_{r_lB}$.
\end{Lm}

\begin{proof}
Put $B^{\rm i}=r_lB$.
Lemma~\ref{goingupequal} applied to $B^{{\rm i}}$ gives that $l \not \sim k$
and that $B^{{\rm i}}$ contains an element $\b$ such that $\b+\a_l$ and
$\b+\a_k$ are in $B$.  Then $B^{\rm i}$ contains both $\b$ and
$\b+\a_l+\a_k$.  This means $(\b,\a_l)=(\b,\a_k)=-1$ in view of
Lemmas~\ref{Prop:AdmOrd} and \ref{heighttworoot}.  Here we distinguish cases
depending on whether or not $j$ is adjacent to $l$ and to $k$.

The easiest case occurs when $j$ is neither adjacent to $l$ nor to $k$.
Here we use the diamond shape from \cite[Lemma $3.1$]{CGW2} with the actions
of $r_l$ and $r_k$.  This gives $B^{\rm iii}=r_jB^{\rm i}<B^{\rm i}$ and
there is a separate path $B>r_jB=B^{\rm ii}>r_lr_jB=B^{\rm iii}$.  As
$r_lB=r_kB$, $B^{\rm iii}=r_kr_jB$ also and $B^{\rm iii}=r_lB^{\rm
ii}=r_kB^{\rm ii}$.  Using induction for the blocks below $B$ we find

\begin{eqnarray*}
r_la_{B^{\rm i}}&=&r_lr_ja_{B^{\rm iii}}\homog r_jr_la_{B^{\rm iii}} \homog r_jr_ka_{B^{\rm iii}}  \homog r_kr_ja_{B^{\rm iii}} = r_ka_{B^{\rm i}}.
\end{eqnarray*}

Next we consider the case where $k\sim j$ and $j\not \sim l$. Here the
following diagram is of use.

$$
\begin{matrix}\label{riBisrjB}
&&B&&&&\\
&l,k\swarrow&&\searrow j&&&\\
B^{\rm i}&&&&B^{\rm ii}&&\\
&j\searrow&&\swarrow l&&\searrow k&\\
&&B^{\rm iv}&&&&B^{\rm iii}\\
&&&k\searrow&&\swarrow j,l&\\
&&&&B^{\rm v}&&\\
\end{matrix}
$$

Set $B^{\rm i}=r_lB=r_kB$, $B^{\rm ii}=r_jB$, $B^{\rm
iii}=r_kB^{\rm ii}$, $B^{\rm iv}=r_jB^{\rm i}$ and $B^{\rm v}=r_kB^{\rm
iv}$.  We use the diamond shape for the actions of $r_l$ and $r_j$ and the
hexagon shape for the actions of $r_k$ and $r_j$ from
\cite[Lemma~$3.1$]{CGW2}.  We use the diamond shape for the actions of $r_l$
and $r_k$ to see $B^{\rm iii} =r_lB^{\rm v}=r_jB^{\rm v}$.  Now we use
induction for the various blocks other than $B$ as they all have lower
height.
\begin{eqnarray*}
r_ka_{B^{\rm i}}&\homog& r_kr_ja_{B^{\rm iv}}\homog r_kr_jr_ka_{B^{\rm
              v}}\homog r_jr_kr_ja_{B^{\rm v}}    \\
              &\homog& r_jr_kr_la_{B^{\rm v}}\homog r_lr_jr_ka_{B^{\rm v}} \homog r_la_{B^{\rm i}}.
\end{eqnarray*}

The final case is $j \sim k$ and $j \sim l$.  We do this much the same as
the above cases but only sketch the argument.  Let $B^{\rm i}=r_kB=r_lB$.
Now let $B^{\rm iv}=r_jB^{\rm i}$.  From here consider the two paths to
$B^{\rm vi}=r_kr_lB^{\rm iv} $ given by $r_l$ and $r_k$.  It is possible
$r_l$ and $r_k$ act the same and this is just one step.  We assume it is
two; the case of just one being easier.  As before we let $B^{\rm ii}=r_jB$.
Again use induction for the blocks other than $B$ which have lower height.
Now
\begin{eqnarray*}
 r_la_{B^{\rm i}}  &\homog& r_lr_jr_lr_ka_{B^{\rm vi}} \homog r_jr_lr_jr_ka_{B^{\rm vi}} \homog   r_ja_{B^{\rm ii}}   \\
              &\homog &       r_jr_kr_jr_la_{B^{\rm vi}}   \homog
              r_kr_jr_kr_la_{B^{\rm vi}} \homog r_ka_{B^{\rm i}}.
\end{eqnarray*}
\end{proof}

This takes care of Case (ii) with $l\not\sim k$, unless there is no $j$
as in Lemma~\ref{equallowerings}.  Assume there is no such $j$.  A search of
all $B\in WB_Y$ for all $Y\in\MY$ using GAP, \cite{GAP}, shows that then
$B$ contains a simple root, say $\a_i$.

We need to show that $r_la_{B^{{\rm i}}} \isog r_ka_{B^{{\rm i}}}$.  As in
the proof of Case (iii) above we may reduce both words all the way
down via sets of the form $B\setminus \{\a_i\}$ for $B\in \B^i$ only
and using elements of $Z_i$
only.  By Lemma~\ref{heightsthesame} they are both
reduced, and as in Case (iii), we find $r_la_{B'}\isog r_ka_{B'}$.

This finishes Case (ii) with $l\not\sim k$.  Next assume $l\sim k$ with
$r_lB<B$ and $r_kB<B$.  Then by Lemma~\ref{Prop:AdmOrd}(iv) either
$r_lr_kB=r_kB$ or the same argument produces paths down starting
$B>r_lB>r_kr_lB<r_lr_kr_lB$.  Using $l$ and $k$ reversed gives an
alternative path through $r_lr_kr_lB$ which can be compared as above.

By Lemma~\ref{heighttworoot}, the case $r_lr_kB=r_kB$ occurs because
$\b=\a_l+\a_k\in B$.  Here $r_l\b=\a_k$ and $r_k(\b)=\a_l$.  An example with
$M=\E_6$ is $Y=\{6\}$ and $\b=\a_5+\a_6$.  If $r_5$ is used
$a_{\{\b\}}=r_5a_{\{\a_6\}}=r_5\he_6$.  If $r_6$ is used
$a_{\{\b\}}=r_6a_{\{\a_5\}}=r_6e_5\he_6$.  We use (HNree) to see
$r_ke_l\isog r_le_ke_l$.  This is sufficient as an alternative to
Definition~\ref{def:aB} of $a_B$ can be made by first taking the product of
$\he_i$ over all nodes $i$ with $\a_i\in \Simp(B)$ and then when a new
simple root $\a_j$ appears in the usual definition after action by $r_l$
multiplying by $\he_j$.  Once there are $|Y|$ different $\he_i$ use the
Temperley--Lieb words as usual.  Then there is no need to multiply by
$\he_Y$ in the final step.

By construction, $a_BB_Y = B$.  If the simple reflection $r_i$ occurs in the
word $a_B $, say $a_B = x r_i y$ for certain words $x$, $y$, then $r_i$
increases the height of the admissible set $xB_Y$ by one. Therefore $\het(B)
=\het(a_B)$.  To finish the proof of Theorem \ref{prop:aB}, observe that
$\pi(\he_Y)\emptyset = B_Y$, so indeed $\pi(a_B)\emptyset = \pi(a_B)B_Y =
B$.

\end{section}

\begin{section}{Reduction to the minimal elements}
\label{sec:aBcharacterization}
This section is devoted to the proof of Theorem \ref{th:aBcharacterization}
for admissible sets $B$.  We use induction and assume the truth of the
theorem for admissible sets of level smaller than $L(B)$ and the truth of
Theorem \ref{prop:aB} for admissible sets of level smaller than or equal to
$L(B)$.

We now begin the proof of Theorem \ref{th:aBcharacterization}.
Let $Y\in\MY$ and $B\in WB_Y$.
We have dealt with the case $\het(B)=0$ in Corollary \ref{cor:aBszero} and
so are assuming that $\het(B) >0$.
Fix a node $i$ of $M$.
We first prove property (i) and next (ii) and (iii) simultaneously.

Here, and later, we will write $=_{\df}$ to indicate that the equality follows
from the definition of $a_B$. Similarly, $\homog_{\rels}$ will indicate that
the reduction is a consequence of the defining relations, and $\homog_{\ih}$
will signify that the reduction is a consequence of the induction hypotheses.

\nl(i).
Recall that, if $B\in WB_Y$, we have $r_iB>B$, $r_iB<B$, or $r_iB=B$.
We treat these cases separately.

If $r_iB>B$, then $a_{r_iB}=_{\df}
r_ia_B$ by the definition of $a_{r_iB}$ and so the result is correct, with
$h$ being the identity, $\he_Y$, of $H_Y$.

If $r_iB<B$, then $ a_B =_{\df} r_ia_{r_iB}$,
so $r_i a_B =  r_ir_ia_{r_iB}\homog_{\rels} a_{r_iB}$, as required with
again $h$ being the identity of $H_Y$.

Suppose then $r_iB=B$.
Now $\a_i$ is perpendicular to all roots of $ B\setminus\{\a_i\}$.
If (iii) of Definition \ref{def:aB} prevails,
there are nodes $j,k$ with $j\sim k$
with $a_B=e_ja_{e_kB}$, $B=e_je_kB$,  and $L(e_kB)<L(B)$.  Notice
$r_iB=B$ implies $\a_i\perp B\setminus\{\a_i\}$.
As $\a_j\in B$, we know $i\not\sim j$.
Now $r_ia_B =_{\df} r_ie_ja_{e_kB} \isog_{\rels}
e_jr_ia_{e_kB}\homog_{\ih} e_ja_{r_ie_kB}h$
for some $h\in H_Y$.

Clearly we are done if $i=j$ using $r_ie_i\homog_{\rels}e_i$ in the first
equality, as then $r_ia_B =_{\df} r_ie_ia_{e_kB} \homog_{\rels} e_ia_{e_kB}
\homog_{\ih} a_{B} h'$ for some $h'\in H_Y$.

Therefore, we may assume $j\ne i$ and (still) $j\not\sim i$.  If $i\not \sim
k$, then $r_ie_kB=e_kB$ and $r_ia_B\homog e_ja_{e_kB}h=_{\df}a_Bh$ and we
are done.

Suppose $i\sim k$.  Notice $\a_k\in e_kB$ and by
Lemma~\ref{sameheight}, $\het(e_kB)=\het(e_ie_kB)$ and
$r_ie_kB> e_kB$ (as $\a_i+\a_k\in r_ie_kB$ and $\a_k\in e_kB$).
We claim $L(e_ie_kB)< L(e_kB)$.  This is because by Definition~\ref{def:aB}
there is a $\b$ in $B$ of minimal height greater than $1$ moved by $r_k$,
for which $(\b,\a_k)=1$ and $\b-\a_k-\a_j\in e_kB$.  Now this is a root of
minimal height moved by $r_i$, is lowered by $r_i$ and so
$L(e_ie_kB)<L(e_kB)$.  We also claim $r_je_ie_kB=e_ie_kB$; for the elements in
$e_ie_kB$ are either perpendicular to $\a_i$, $\a_j$, and $\a_k$ or of the
form $\c+2\eps\a_k +\eps \a_j +\eps \a_i$ where $\c\in B$ and
$\eps=-(\c,\a_k)$. (To see this, use the action of $e_k$ on $B$ to be
$r_jr_k$ and the action of $e_i$ on $e_B$ to be $r_kr_i$).  Now $\a_j$ is
orthogonal to these.  Notice also that $r_k(e_ie_kB)>e_ie_kB$, as the root
$\b-2\a_k -\a_j - \a_i$, for $\b$ as above, is a root of
$e_ie_kB$ of minimal height moved by $r_k$ and is raised by $r_k$.
Now we have enough properties to conclude

\begin{eqnarray*}
e_ja_{r_ie_kB}&=& e_ja_{r_ke_ie_kB}
\isog_{\ih}    e_jr_ka_{e_ie_kB}
\homog_{\rels}  e_je_kr_ja_{e_ie_kB}\\
&\homog_{\ih} &
 e_je_ka_{e_ie_kB}h'\homog_{\ih}
 e_ja_{e_ke_ie_kB}h'' = e_ja_{e_kB}h''\\
&\homog_{\df} &
a_{B}h'',
\end{eqnarray*}
for certain $h',h''\in H_Y$ and so
$r_i a_B\homog a_Bhh''=a_{r_iB}hh''$, as required.
This settles the case where Definition \ref{def:aB}(iii) applies.

Suppose next (ii) of Definition \ref{def:aB} prevails, that is, there is a
node $k$ of $M$ such that
$a_B=_{\df}r_ka_{r_kB}$ with $\het(r_kB)<\het(B)$.
We know $i\neq k$ as $r_iB=B$.

Assume $i\not\sim k$. Now $r_ir_kB = r_kr_iB = r_kB$, so there is $h\in H_Y$
such that $r_i a_{B} =_{\df} r_ir_ka_{r_kB} \isog_{\rels} r_k r_i
a_{r_kB}\homog_{\ih} r_k a_{r_kB}h =_{\df} a_{B}h $, as required.

Assume $i\sim k$. Then $r_ir_kB = r_ir_kr_i B = r_k(r_ir_k B)$, so $r_k$
fixes $r_ir_kB$.  By definition $r_kB<B$ and so $r_k$ raises $r_kB$.
This
means that $r_k$ raises all of the elements in $r_kB$ of smallest height
that are moved by $r_k$.  Such a root $\b\in r_kB $ is moved to $\b+\a_k\in B$ under
the action of $r_k$.
As $r_iB=B$, we have $r_i(\b+\a_k)=\b+{\a_k}$ and so $(\b,\a_i)=1$.  This
means $r_i$ lowers the elements of smallest height of $r_kB$ that $r_k$
raises.  Elements of $r_kB$ not moved by $r_k$ are not moved by $r_i$ and so
$r_i$ lowers $r_kB$ and we can use induction.  This gives $h\in H_Y$ such
that $r_i a_{B} =_{\df} r_i r_k a_{r_kB}\isog_{\ih} r_i r_k r_i a_{r_i
r_kB} \homog_{\rels} r_kr_i (r_k a_{r_ir_kB}) \homog_{\ih} r_kr_i
(a_{r_ir_kB}) h \homog_{\ih} r_k (a_{r_kB}) h \homog_{\df} a_{B}h =
a_{r_iB}h $, as required.

We have dealt with cases (ii) and (iii) of Definition \ref{def:aB}.  In case
(i), the height of $B$ is zero, so by our assumption $\het(B)>0$, all
possibilities are exhausted and the induction step for
Theorem \ref{th:aBcharacterization}(i) is proved.

\medskip
We now come to the proof of the induction step for (ii) and (iii) of Theorem
\ref{th:aBcharacterization}.  We will deal with these
simultaneously, proceeding in a number of steps.
By using GAP, we are able to show that all cases are eliminated proving the
theorem as we describe at the end of this section.

\begin{Remark} \rm
In many instances we have $\het(B)=\het(e_iB)$.  In these
cases (ii) can be improved to $e_ia_B \isog a_{e_iB}$ with no $h$ appearing.
This is because both $a_B$ and $a_{e_iB}$ have the same height, $\het(B)$,
and are both reduced.  This means that $h$ is the identity $\he_Y$ of $H_Y$.
We use this sometimes without referring to it.
\end{Remark}

The first several of these steps concern the case where $j$ is a lowering
node for $B$, so $r_jB<B$.  This implies $a_B = r_ja_{r_jB}$.  Notice that
by the induction assumptions any two definitions for $a_B$ must be the same
up to $\isog$ as each will be reduced of the same height, $\het(B)$.

\begin{Step}\label{test1}
Suppose $r_jB<B$, and
$j\not\sim i$.  Then
(ii) and (iii) hold.
\end{Step}

\begin{proof}
Here $a_B=_{\df} r_ja_{r_jB}$ by definition and $e_ir_ja_{r_jB}\isog_{\rels}
r_je_ia_{r_jB}$.  If $|e_ir_jB|> |r_jB|$, then by induction, the word
$e_ia_{r_jB}$ reduces to an element as in (iii) and so does $e_ia_B$.
Therefore, we may assume $|e_ir_jB|=|r_jB|$.
Then, again by induction, we find
$\het(e_ir_jB)\le \het(r_jB)<\het(B)$ and
there are $h,h'\in H_Y$ such that
$r_je_ia_{r_jB}\homog_{\ih} r_ja_{e_ir_jB}h\homog_{\ih}
a_{r_je_ir_jB}h' =
a_{e_iB}h'$, so $e_ia_B\homog a_{e_iB}h'$.

\end{proof}

In the remaining steps these checks for $\het (e_iB)\leq \het(B)$ when
$|e_iB|=|B|$ are routine and we leave them to the reader.
With the exception of Step~\ref{test13}, we do the same when in a step in
the induction we have an instance of $|e_jB'|>|B'|$ for a $B'$ of lower
height with an $e_j$ appearing in a step, leading to an instance of (iii).

Also, often reduction steps are written down without the explicit mention of
powers of $\delta$ that might occur as factors. They are dropped for the
sake of simplicity as they have no bearing on the result.

\begin{Step}\label{test8} Suppose $r_jB<B$ and $r_ir_jB<r_jB$.
Then  (ii) and (iii) hold.
\end{Step}

\begin{proof}
In view of Step \ref{test1}, we may assume $i\sim j$.  Notice
$e_ia_B\isog_{\ih}e_ir_jr_ia_{r_ir_jB}\homog_{\rels}
e_ie_ja_{r_ir_jB}$. (The absence of elements from $H_Y$ is due to the second
statement of Theorem \ref{th:aBcharacterization}(i).) As $e_je_iB =
e_j(r_ir_jB)$, we have, by induction $\het(e_je_iB)<\het(B)$. Now use
induction to find $h,h'\in H_Y$ with
\begin{eqnarray*} e_ie_ja_{r_ir_jB}&\homog_{\ih}
  &e_ia_{e_jr_ir_jB}h=e_ia_{e_je_iB}h\homog_{\ih} a_{e_ie_je_iB}h' =  a_{e_iB}h',
\end{eqnarray*}
 as required.  As mentioned, we are leaving to the reader the cases in which $|e_jr_ir_jB|>|r_ir_jB|$ and $|e_ie_je_iB|>|e_je_iB|. $
\end{proof}

\begin{Step}\label{test10}  Suppose $r_jB<B$ and $r_ir_jB=r_jB$. Then
(ii) and
(iii) hold.
\end{Step}

\begin{proof}
The case $i\not\sim j$ is dealt with by Step \ref{test1}, so without loss of
generality, we assume $i\sim j$.  Using the definition, the relations, and
induction $e_ia_B=_{\df}e_ir_ja_{r_jB}\isog_{\rels}
e_ie_jr_ia_{r_jB}\homog_{\ih} e_ie_ja_{r_ir_jB}h=e_ie_ja_{r_jB}h$
for some $h\in H_Y$.
Now use the
induction twice to find $h',h''\in H_Y$ with $e_ie_ja_{r_jB}\homog_{\ih}
e_ia_{e_jr_jB}h'\homog_{\ih} a_{e_ie_jB}h''$.  Now,
$e_ie_jB=e_ir_jr_iB = e_ir_ir_jr_ir_jB = e_ir_jr_jB =  e_iB$,
so $ e_ie_jB=e_iB$ and we
are done.
\end{proof}

\begin{Step}\label{test14}  Suppose that $r_jB<B$ and
\begin{itemize}
\item[(a)] $e_ja_{r_ir_jB}\homog a_{e_jr_ir_jB}h'$ and
\item[(b)] $e_ia_{e_jr_ir_jB}\homog a_{e_ie_jr_ir_jB}h''$
\end{itemize}
both hold for $h', h'' \in H_Y$.
Then (ii) and
(iii) hold.

A sufficient condition for (a) to hold is that there is a node $k$ with
$k\not\sim j\sim i$ such that $r_k$ lowers $r_ir_jB$. A sufficient condition
for (b) to hold is that there is a node $l$ with $l\not \sim i$ that lowers
$e_jr_ir_jB$ or that $L(e_jr_ir_jB) < L(B)$.
\end{Step}

\begin{proof}
As for the first assertion, in view of Step~\ref{test1} and the definition
we may assume $i\sim j$.  Using part (i) and induction we see
$e_ia_B=_{\df}e_ir_ja_{r_jB}\homog_{\rels} e_ie_jr_ia_{r_jB}\homog_{\ih}
e_ie_ja_{r_ir_jB}$.  Because (a) and (b) both hold, this reduces to
$a_{e_ie_jr_ir_jB}h''h'$.  As $e_ie_jr_ir_jB = e_ie_je_iB = e_iB$ by (RNerr)
and (HNeee), the result follows.

As for the second assertion, the hypothesis on $k$ implies
$e_ja_{r_ir_jB}\homog a_{e_jr_ir_jB}h'$ for some $h'\in H_Y$ by
Step~\ref{test1}, which means (a) holds.

As for the conditions for (b), the condition $L(e_jr_ir_jB) < L(B)$
implies (b) by induction.  If $l$ lowers $e_jr_ir_jB$ and $i\not\sim l$,
then $e_ia_{e_jr_ir_jB}\homog a_{e_ie_jr_ir_jB}h''$
for some $h''\in H_Y$ also by Step~\ref{test1},
which means (b) holds.  This finishes the step.
\end{proof}

In the next three steps there may or may not be a lowering node for $B$.

\begin{Step}
\label{test7} Suppose there are no lowering nodes for $e_iB$ and
$k$ is a lowering-e-node for $e_iB$ with
$\a_k\in B$ and $k\sim i$.  Then $e_ia_B\isog a_{e_iB}$.
\end{Step}

\begin{proof}  By Lemma~\ref{sameheight}, we have
$e_ke_iB=B$ and Definition \ref{def:aB}(iii) with $L(e_ke_iB)<L(e_iB)$ gives
$a_{e_iB}=_{\df}e_ia_{e_ke_iB} = e_ia_{B}$, as required.
\end{proof}


\begin{Step}\label{test9}  Suppose there are no lowering nodes for $e_iB$
and $j$ is a node with $\a_j\in B$ and $i\sim j$.
Suppose also $L(e_iB)<L(B)$.  If either
there is a node $k$ with
$L(e_ke_iB)<L(e_iB)$ and $i\sim k$, or
$B$ has no lowering
nodes,  then $e_ia_B\isog a_{e_iB}$.
\end{Step}

\begin{proof}
Suppose first there is a node $k$ as indicated.  Using the definition, the
relations, and induction we see $a_{e_iB}=_{\df}e_ia_{e_ke_iB}\isog_{\rels}
e_ie_je_ia_{e_ke_iB} \homog_{\ih} e_ie_ja_{e_ie_ke_iB}=e_ie_ja_{e_iB}$.
Notice there is no $h$ term here as ${e_iB}$ and $B $ are of the same
height in the poset, and $a_{e_iB}$ and $a_B $ are reduced of this same
height. This means that we even have
$a_{e_iB}\isog e_ie_ja_{e_iB}$.  Now use induction to see
$e_ie_ja_{e_iB}\homog_{\ih} e_ia_{e_je_iB}=e_ia_{B}$. By the same argument
as before, we may replace the occurrence of $\homog $ by $\isog$, and so we
are done.

Suppose now $B$ has no lowering nodes.  Then $a_B=_{\df}e_ja_{e_iB}$.
Now $e_ia_B=e_ie_ja_{e_iB}\homog_{\ih}
e_ie_je_ia_{e_iB}\homog_{\rels}e_ia_{e_iB}\homog_{\ih}a_{e_iB}$. Again the
occurrences of $\homog$ can be
replaced by $\isog$, which leads to the required result.
\end{proof}

\begin{Step}\label{test11}
Suppose $L(e_jB) < L(B)$.
If $k$ is a node with $\a_k \in B$ satisfying
$i\not\sim k \sim j$
and $L(e_ie_jB) < L(B)$,
then
(ii) and (iii) hold.
\end{Step}

\begin{proof}
We have $B = e_ke_jB$ and $\het(e_ke_jB) = \het(e_jB)$,
so there are $h,h'\in H_Y$ with
\begin{eqnarray*}
 e_i a_B &=& e_i a_{e_ke_jB} \isog_{\ih} e_i e_k a_{e_jB} \isog
_{\rels}e_k e_ia _{e_jB} \\
&\homog_{\ih}& e_k a _{e_ie_jB}h \homog_{\ih}  a _{e_ke_ie_jB}h' =
   a _{e_ie_ke_jB}h'  \\
&=&   a _{e_iB}h',
\end{eqnarray*}
 as required.
\end{proof}

For the remainder of the proof we may assume there is no node $j$ with
$r_jB<B$.  This means that Definition \ref{def:aB}(iii) applies and there
are adjacent nodes $j$, $k$ with $\a_k\in B$ and $a_B = e_k a_{e_jB}$.

\begin{Step}\label{test2} Suppose  $j$ is a lowering-e-node of $B$ with
$i\sim j$.  If $\a_i\in B$, then (ii) and (iii) hold.
\end{Step}

\begin{proof}
By Lemma \ref{sameheight}, $e_ie_jB = B$ and  $B=e_iB$.
By definition $e_ia_B=_{\df}e_ie_ia_{e_jB}$.
As $L(e_jB)<L(B)$ we can use induction and, as $e_i^2\isog_{\rels}\delta
e_i$, we find $h\in H_Y$ with $e_ie_ia_{e_jB}\isog_{\rels}\delta
e_ia_{e_jB}\homog_{\ih}
a_{e_ie_jB}h =\delta a_B h=\delta a_{e_iB}h$.
\end{proof}

\begin{Step}\label{test3} Suppose  $i$ is a lowering-e-node for $B$.  Suppose
$j\sim i$ with $e_je_iB=B$ and $L(e_iB)<L(B)$.  Suppose also
$k$ lowers $e_iB$ and $i\not \sim k$.  Then $e_ia_B\isog a_{e_iB}$.
In particular, (ii) and (iii) hold.
\end{Step}

\begin{proof}
Using the definition and induction, we find
$e_ia_B=_{\df}e_ie_ja_{e_iB}\isog_{\df}e_ie_jr_ka_{r_ke_iB}$.  As
$e_i^2B=e_iB$, we find
\begin{eqnarray*}
e_ie_jr_ka_{r_ke_iB}&=&
e_ie_jr_ka_{r_ke_i^2B}=_{\rels} e_ie_jr_ka_{e_ir_ke_iB}\isog_{\ih}
e_ie_jr_ke_ia_{r_ke_iB}\\
&\isog_{\rels}& e_ie_je_ir_ka_{r_ke_iB}=_{\rels}
e_ir_ka_{r_ke_iB}\isog_{\ih}e_ia_{e_iB}\\
&\homog_{\ih}&a_{e_iB},
\end{eqnarray*}
where the
absence of factors $h\in H_Y$ is explained as before and the last induction
step is valid because $L(e_iB)<L(B)$.
\end{proof}

Recall $N(j,B)$ from Notation \ref{not:NlB}.

\begin{Step}\label{test4} Suppose $j$ is a lowering-e-node for $B$
and $k=N(j,B)$.
If
$k\sim i\not \sim j$, and $i$ is a lowering node for $e_jB$,  then
(ii) and (iii) hold.
\end{Step}

\begin{proof}
Using Definition \ref{def:aB}(iii), (ii),
we see $e_ia_B=_{\df}e_ie_ka_{e_jB}=_{\df}e_ie_kr_ia_{r_ie_jB}
\homog_{\rels} e_ir_ka_{r_ie_jB}$.
Notice that $r_kr_ie_jB$ has the same height as $B$ as $k$ raises the
simple root $\a_j$ in $r_ie_jB$ to $\a_j+\a_k$ in $r_kr_ie_jB$.  Now $r_j$
moves $\a_j+\a_k$ to $\a_k$ and so
$r_jr_kr_ie_jB<r_kr_ie_jB$. Therefore, there are $h,h',h''\in H_Y$ such that
\begin{eqnarray*}
e_ir_ka_{r_ie_jB}  &\homog_{\ih} & e_ia_{r_kr_ie_jB}h \homog_{\df} e_ir_ja_{r_jr_kr_ie_jB}h
        \homog_{\rels}       r_je_ia_{r_jr_kr_ie_jB}h  \\
          &\homog_{\ih}&  r_ja_{e_ir_jr_kr_ie_jB}h'
      \homog_{\ih}  a_{r_je_ir_jr_kr_ie_jB}h''.
      \end{eqnarray*}
Now $r_je_ir_jr_kr_ie_jB=e_ir_kr_ie_jB=e_ie_ke_jB=e_iB$, which finishes the proof.
\end{proof}

\begin{Step}\label{test13}   Suppose that $k$ is a lowering-e-node for $B$
and $j=N(k,B)$ satisfies $j\not\sim i$.
If $|e_ie_kB|>|B|$,  then $e_ia_B$ reduces to an element of
$\Br(M)e_U\Br(M)$ for some $U$ properly containing $Y$, so
(ii) and (iii) hold.
\end{Step}

\begin{proof}
Notice $e_ia_B =_{\df} e_ie_ja_{e_kB}\homog_{\rels}e_je_ia_{e_kB}$.  Now as
$L(e_kB)<L(B)$, induction together with $|e_ie_kB|>|B|$ gives that
$e_ia_{e_kB}$ reduces to an element as stated, and hence
$e_je_ia_{e_kB}$ as well.
\end{proof}

All possible instances of reduction of $e_ia_B$ as in (ii) and (iii) for
$M\in\{\E_6,\E_7,\E_8\}$ are covered by Steps~\ref{test1}
to \ref{test13}.  This fact has been checked by use of GAP \cite{GAP}.

\end{section}

\begin{section}{Conclusion}\label{sec:conclusion}
In this section we prove Theorem~\ref{th:main}.  To this end, we establish
Theorem~\ref{th:nottheta} (in \ref{proof:nottheta}) as a consequence of the
results in the previous sections. Then we derive part (i) of Theorem
\ref{th:main}.  Next we will be concerned with semisimplicity (Theorem
\ref{th:ss}) and cellularity (Theorem \ref{th:cellular}), proving the
remaining parts, (ii) and (iii), of the same theorem.  These two properties
are established in much the same way the corresponding result is shown for
$\D_n$ in \cite[Section~5]{CGW3}.  We conclude with a remark on subalgebras
associated with subdiagrams of $M$.

\begin{Proof} {\bf Proof of Theorem \ref{th:nottheta}.}
\rm
\label{proof:nottheta}
As before, it suffices to deal with the cases $M=\E_n$ $(n=6,7,8)$.  Suppose
$a\in F$ and write $B = \pi(a)(\emptyset)$ and $B' =
\pi(a^{\op})(\emptyset)$. Let $Y\in\MY$ be such that $B\in WB_Y$.  We need
to show that $a$ can be reduced to an element of the form $\delta^i a_B\he_Y
h \a^{{\op}}_{B'}$ for some $i\in\Z$ and $h\in T_Y$.  The existence of $a_B$
is established in Theorem \ref{prop:aB} and is unique $\tF$ up to powers of
$\delta$ by the same theorem.

We do so by induction on the length of $a$ in terms of the generators $r_i$
and $e_i$ (and so disregarding the powers of $\delta$). If $a$ is the empty
word $1$, then clearly $a(\emptyset) =a^{\op}\emptyset = \emptyset$, and $a
= a_{\emptyset}1a_{\emptyset}^{\op}$, so the theorem holds.

Now suppose $a = x b$ with $x$ a generator of $F$ different from
$\delta^{\pm1}$ and $b$ a word of $F$. Then, by induction on the length of
$a$, there are a subset $Y'$ of $Y$, admissible sets $C,C'\in WB_{Y'}$ and
$y\in T_{Y'}$ such that $b\homog\delta^j a_{C}ya_{C'}^{\op}$. If $x = r_i$
then $B = r_iC$ and $B'=C'$, so $Y'=Y$ and, by Theorem
\ref{th:aBcharacterization}(i), there is $h\in H_Y$ such that $a \homog
\delta^j a_{B} hya_{C'}^{\op}$, and we can finish by Theorem \ref{th:WC},
which gives us that we may in fact assume $hy\in T_Y$.

Next suppose $x = e_i$ for some node $i$ of $M$. If $e_iC\in WC$, then we
can argue as for $x=r_i$, using Theorem \ref{th:aBcharacterization}(ii).
So, we may assume $\a_i\perp C$ and $B = e_iC = \cl{(C\cup\{\a_i\})} $. Now,
Theorem \ref{th:aBcharacterization}(iii) and repeated application of the
other parts of the theorem give $v\in F$ such that $a =e_ib \homog
\delta^{j}a_{B}\he_Y v ya_{C'}^{\op}$ with $B_Y\pi(vya_{C'}^{\op}) = B'$ for
some $j\in\Z$. Another application of parts (i) and (ii) of Theorem
\ref{th:aBcharacterization} and of Theorem \ref{th:WC} gives $w\in T_Y$ such
that $ a_{C'}(vy)^{\op}\he_Y \homog \delta^{k}a_{B'} w $ for some $k\in\Z$.
We conclude $a\homog \delta^j a_{B}\he_Y v ya_{C'}^{\op} \homog\delta^{j+k} a_{B}\he_Y w^{\op}
a_{B'}^{\op} \homog \delta^{j+k} a_{B} w^{\op} a_{B'}^{\op}$, as required.
\end{Proof}

\begin{Proof} {\bf Proof of Theorem \ref{th:main}(i)}.
\label{proofMainTheorem}
\rm Choose a set $T$ of words in $F$ whose image under $\pi$ is a set of
representatives for the regular group action of $\langle\delta\rangle$ on
$T_\delta$, as described in (\ref{dfT}).  Then, by Theorem
\ref{th:nottheta}, each word in $\tF$ reduces to a unique element of $T$ up
to a power of $\delta$.
By Proposition~\ref{prop:BMWbasis} the
set $\rho(T)$ is a basis for $\BMW(\E_n)$ and so $\BMW(\E_n)$ is free of the
correct rank.  This
proves Theorem \ref{th:main}(i).
\end{Proof}

\begin{Thm}\label{th:ss}
If $M=\E_n $ for $n\in \{6,7,8\}$, then  $\BMW(M)\otimes_R
\Q(l,\delta)$ is semisimple.
\end{Thm}

\begin{proof}
To show that $\BMW(M)$ tensored over $\Q(l,\delta)$ is semisimple we use
the surjective ring homomorphism $\mu\,\colon\BMW(M)\otimes _R
\Q(\delta)[l^{\pm1}]\to\Br(M)$ over $\Q(\delta)$ defined in
Section~\ref{intro} just after the proof of Proposition~\ref{templieb}.  We
know its image $\Br(M)$ is semisimple by \cite[Corollary~5.6]{CFW} and so
has no nilpotent left ideals.  Suppose $\BMW(M)\otimes _R\Q(\delta,l)$
has a nontrivial nilpotent ideal.  Take a nonzero element of it expressed in
the basis we have found.  Multiply the element by a suitable polynomial in
$l$ so that all coefficients are in $\Q(\delta)[l^{\pm 1}]$.  As in the
proof of \cite[Lemma~4.2]{CGW3}, rescale the coefficients by a power of
$l-1$ so that all coefficients remain in $\Q(\delta)[l^{\pm1}]$ but some
coefficient $\lambda_s$ lies outside $(l-1)\Q(\delta)[l^{\pm1}]$.  The
result is a nonzero nilpotent element in $\BMW(M)\otimes
\Q(\delta)[l^{\pm1}]$ with $\mu(\lambda_s) \neq 0$, so its image under $\pi$
is a nonzero nilpotent element of $\Br(M)$.  Furthermore, any multiple is
nilpotent both in $\BMW(M)\otimes \Q(\delta,l)$ and in $\Br(M)$ and so
generates a nontrivial nilpotent ideal of $\Br(M)$, a contradiction with
the semisimplicity of $\Br(M)$.  This completes the proof of Theorem
\ref{th:main}(ii).
\end{proof}

\begin{Remark} \rm
By use of $\mu$ and the Tits Deformation Theorem, see \cite[IV.2, exercice
26]{Bour} or \cite[Lemma 85]{Stein}, it can be shown that the irreducible
degrees associated to $\BMW(\E_n)$ are the same as for $\Br(\E_n)$ for $n=6,7,8$.
\end{Remark}

Next we prove the cellularity part of Theorem~\ref{th:main}.
The proof given here runs in the same way as the proof of the
corresponding result for $\D_n$ in \cite[Section~6]{CGW3}.
The result is stated in Theorem \ref{th:cellular}.

Recall from \cite{GL} that an associative algebra $\alg$ over a commutative
ring $S$ is cellular if there is a quadruple $(\Lambda, D, C, *)$ satisfying
the following three conditions.

\begin{itemize}
\item[(C1)] $\Lambda$ is a finite partially ordered set.  Associated to each
$\lambda \in \Lambda$, there is a finite set $D(\lambda)$.  Also, $C$ is an
injective map
$$ \coprod_{\lambda\in \Lambda} D(\lambda)\times D(\lambda) \rightarrow \alg$$
whose image is an $S$-basis of $\alg$.

\item[(C2)]
The map $*:\alg\rightarrow \alg$ is an
$S$-linear anti-involution such that
$C(x,y)^*=C(y,x)$ whenever $x,y\in
D(\lambda)$ for some $\lambda\in \Lambda$.

\item[(C3)] If $\lambda \in \Lambda$ and $x,y\in D(\lambda)$, then, for any
element $a\in \alg$,
$$aC(x,y) \equiv \sum_{u\in D(\lambda)} r_a(u,x)C(u,y) \
\ \ {\rm mod} \ \alg_{<\lambda},$$ where $r_a(u,x)\in S$ is independent of $y$
and where $\alg_{<\lambda}$ is the $S$-submodule of $\alg$ spanned by $\{
C(x',y')\mid x',y'\in D(\mu)\mbox{ for } \mu <\lambda\}$.
\end{itemize}
Such a quadruple $(\Lambda, D, C, *)$ is called a {\em cell datum} for
$\alg$.
We will describe such a quadruple.  For $*$ we will use $^\op$ defined by
\begin{Notation}\label{not:DX}
\rm
For
$x_1,\ldots, x_q\in\{r_1,\ldots,r_n,e_1,\ldots,e_n,\delta^{\pm1}\}$, we write
$(x_1\cdots x_q)^\op = x_q\cdots x_1$, thus defining an opposition map on
$F$. This notation is compatible with the maps $\pi$ and $\rho$ when
${\cdot}^\op$ on $\BMW(\E_n)$ and $\Br(\E_n)$ is interpreted as the
anti-involution of \cite{CGW}
and \cite{CFW}, respectively; see Definition \ref{def:Baction}.
\end{Notation}

\np We introduce a quadruple $(\Lambda, D,C,*)$ and prove that it is a cell
datum for $\alg = \BMW(\D_n)\otimes_R S$.  Before describing these, we will
relate the subalgebras of $\alg$ generated by monomials corresponding to the
elements of $S_Y$ in Table~\ref{table:types2} to Hecke algebras.  Here for
$S_Y$ in Table~\ref{table:types2} we act on each term by $\rho $ to get
elements in $\BMW(\E_n)$. For this purpose we need a version of
Corollary~\ref{Cor:WC} that applies to $\alg$ rather than $\BrM(\E_n)$. This
requires a version of Theorem~\ref{th:WC} for $\BMW(\E_n)$ rather than
$\tF$.  What we do here corresponds to \cite[Corollary~$6.4$]{CGW3}.

In particular we use the following ideals in $\BMW(\E_n)$.

\begin{Def}\label{BrauerJt} \rm
For $M=\E_n$ $(n=6,7,8)$ and a $Y\in\MY$, let $t=|B_{Y}|$ be as listed in
column~$2$ of the row for $Y$ in Table~\ref{table:types2} if $Y\ne\emptyset$
and $t=0$ otherwise.  Put $J_0 = \Br(M)$.  If $t>0$, we let $J_t$ be the
ideal of $\Br(M)$ generated by $e_Y$ together with all $e_{Y'}$ for sets of
nodes $Y'$ with $|Y'|>t$.
\end{Def}

In this section we use the same notation for the corresponding ideals in
$\BMW(\E_n)$ and trust it will not cause confusion.

\begin{Def}\label{BMWJt} \rm
For a fixed $Y$ in column~$5$ of Table~\ref{table:types2} we let $J_t$ be the
ideal of $\BMW(\E_n)$ generated by $\rho(e_Y)$ together with all
$\rho(e_{Y'})$ for all $Y'$ in a row lower than $Y$ for that $\E_n$.
Here $t=|B_Y|$
is listed in column~$2$ of the row.
\end{Def}

\begin{Def}\label{def:si}  \rm
Let $Y\in \MY$. For each of the elements of $S_Y$
(see Table~\ref{table:types2}, column~$6$) of the form $r_i\he_Y$, we
let $\hat{s}_i$ be the image $\rho(r_i\he_Y)$ in $\BMW(\E_n)$.
For each of the
remaining long words in $S_Y$ (occurring in column 6 of
rows $1,2,4,5,9,10$) we let
$\hat{s}_0$ be the image of $\rho$ on the word.  In particular
for the first row $\hat{s}_0=e_6e_5e_4g_2e_3e_4e_5e_6\delta^{-1}$.
\end{Def}

\begin{Prop}\label{hecke} For each row of Table~\ref{table:types2}
the $\hat{s}_i$ of Definition~\ref{def:si} are generators of the Hecke algebra of type
$M_{Y}$ which appears in the fourth column of the row modulo the ideal
$J_{t+1}$.  Here $t=|B_Y|$ in the second column of the row.
\end{Prop}

\begin{proof}  To show that the generators $\hat{s}_i$ generate the Hecke algebra of type $M_{Y}$ we need to
show they satisfy the braid relations as well as the quadratic relations.
The proof of this is very much like the corresponding proof in
\cite[Lemma~$6.1$]{CGW3}.  The braid relations for the $S_Y$ have been
proved in Theorem~\ref{th:WC} for $\isog$ which implies they are true in
$\BrM(\E_n)$.  To show they are satisfied in $\BMW(\E_n)$ we have to show
they are still true when the remaining terms involving $m$ occur.  Many of
the relations are binomial terms with no $m$ appearing in
Table~\ref{BMWTable}.  These are all except (RSrr), (HNrer), (HNree), and
(HNeer).

We start with the quadratic terms.  For this we must show
$\hat{s}_i^2=1-m{\hat s}_i$ mod $J_{t+1}$.  Because of (RSrr) for $i\neq 0$
we need to show $ml^{-1}e_i$ acts as $0$.  In these cases $e_ie_Y$ is in
$J_{t+1}$ as $i$ is not adjacent to a node in $Y$.  The other case is
$\hat{s}_0$.  For this we do the case $\hat{s}_0$ for $\E_6$ with
$Y=\{\a_6\}$, so $t=1$.  The main part of this needs
\begin{eqnarray*}
                             \hat{s}_0^2                           &=&  e_6e_5e_4g_2e_3e_4e_5e_6e_6e_5e_4g_2e_3e_4e_5e_6\delta^{-2}  \\
                                                             &\isog& e_6e_5e_4g_2e_3e_4e_5e_6e_5e_4g_2e_3e_4e_5e_6\delta^{-1} \\
                                                             &\isog& e_6e_5e_4g_2e_3e_4e_3g_2e_4e_5e_6 \delta^{-1}               \\
                                                             &\isog&  e_6e_5e_4g_2e_3g_2e_4e_5e_6 \delta^{-1}                      \\
                                                              &\isog&  e_6e_5e_4g_2^2e_3e_4e_5e_6 \delta^{-1}               \\
                                                              &\isog& e_6e_5e_4(1-mg_2+ml^{-1}e_2)e_3e_4e_5e_6\delta^{-1}  \\
                                                               &\isog& e_6e_5e_4e_3e_4e_5e_6\delta^{-1} -m e_6e_5e_4e_3g_2e_4e_5e_6\delta^{-1} +ml^{-1}e_6e_5e_4e_3e_2e_4e_5e_6\delta^{-1} \\
                                                             &\isog &
                                                             1-m\hat{s}_0   \
                                                             \ {\rm mod}\ J_{2}
\end{eqnarray*}
as $e_6e_5e_4e_3e_2e_4e_5e_6$ is in $J_2$ as
$e_3e_2e_4e_5e_6\{\a_6\}=\{\a_3,\a_6\}$.  The braid relations for the elements
not including $s_0$ follow from the ordinary braid relations.  For the ones
containing $s_0$ we have to modify the proof of Theorem~\ref{th:WC} by
including the terms involving $m$.

We do first $\hat{s_2}\hat{s_0}\isog \hat{s_0}\hat{s_2}$ for the case $\E_6$ row~$1$ with $\hat{s_0}=e_6e_5e_4e_3g_2e_4e_5e_6\delta^{-1}$.  This
is covered by Lemma~$7.1$ of \cite{CGW3} but we include the details here with the current notation.
\begin{eqnarray*}
\hat{s_3}\hat{s_0}\delta &\isog & g_3e_6e_5e_4g_2e_3e_4e_5e_6    \\
                         &\isog & e_6e_5g_3e_4e_3g_2e_4e_5e_6        \\
                          &\isog & e_6e_5g_3^2g_4e_3g_2e_4e_5e_6     \\
                          & \isog &  e_6e_5g_4e_3g_2e_4e_5e_6 -me_6e_5g_3g_4e_3g_2e_4e_5e_6  \\
                                        & &\qquad\qquad\qquad  +ml^{-1}e_6e_5e_3g_4e_3g_2e_4e_5e_6  \\
                          &\isog & e_6e_5e_4g_5^{-1}e_3g_2e_4e_5e_6 -me_6e_5e_4e_3g_2e_4e_5e_6    \\
                                       & & \qquad\qquad\qquad +me_6e_5e_3g_2e_4e_5e_6             \\
                          &\isog & e_6e_5e_4g_5^{-1}e_3g_2e_4e_5e_6 -me_6e_5e_4e_3g_2e_4e_5e_6    \\
                                       & & \qquad\qquad\qquad +me_3g_2e_6e_5e_4e_5e_6            \\
                          &\isog & e_6e_5e_4g_5^{-1}e_3g_2e_4e_5e_6 -me_6e_5e_4e_3g_2e_4e_5e_6    \\
                                       & & \qquad\qquad\qquad +me_3g_2e_6
\end{eqnarray*}

Notice that all terms in the last line are fixed under $^\op$ and so $\hat{s}_3\hat{s}_0$ is also and so $\hat{s}_3$ and $\hat{s}_0$ commute.

The other commuting cases in this example are also covered by \cite[Lemma~$7.1$]{CGW3}.

We now tackle the case $\hat{s_1}\hat{s_0}\hat{s_1}\isog
\hat{s_0}\hat{s_1}\hat{s}_0$.  This can be done by the same methods of
computations but the details are messy.  We present another method which
relies on the isomorphism of the BMW algebras of type $\A_{n-1}$ with
tangles on $n$ strands as shown in \cite{MorWas}.

The case we present is really the case for $M = \E_6$ with $|X|=2$.  Here
$\hat{s}_0=e_4e_3g_2e_4e_6\delta^{-2}$ and $\hat{s}_1=g_1e_4e_6\delta^{-2}$.
We do a computation with tangles for $g_1$ and $e_4e_3g_2e_4$ and note this
is sufficient for all of the cases with $\hat{s}_0$ appearing by using
computations which do not introduce extra terms involving $m$.

In particular we show $g_1e_4e_3g_2e_4g_1\isog \delta
e_4e_3g_2e_4g_1e_4e_3g_2e_4$.  After putting in the appropriate $\delta$s
this is what is needed to show $\hat{s_0}\hat{s}_1\hat{s}_0\isog
\hat{s_1}\hat{s}_0\hat{s}_1$.

Notice these elements are all in an $\A_4$ with generators $g_1,g_3,g_4,g_2$
and $e_1,e_3,e_4,e_2$ taken in this order as this order generates an $\A_4$ in
terms of the nodes of the Dynkin diagram we are using.  The tangles then are
on $5$ strands.  For our purposes we take five nodes at the top labelled
$1,3,4,2,5$ arranged horizontally in that order and on the bottom five more
nodes labelled $\bar1,\bar3,\bar4,\bar2,\bar5$ also arranged horizontally in
that order with $i$ directly above $\bar i$ for $i=1,3,4,2,5$.  The tangle for
$g_1$ has $3$ joined to $\bar1 $ and $1$ joined to $\bar 3$ with the strand
from $3$ to $\bar1 $ above the strand from $1$ to $\bar 3$.  The remaining
strands are vertical strands from $i$ to $\bar i$ for $i=4,2,5$.  The tangle
for $e_4$ has $4$ and $2$ joined as well as $\bar4$ and $\bar2$ plus vertical
strands for the remaining vertices $1,3,5$.  The tangle for $e_3$ is similar
except $3$ and $4$ are joined and as well as $\bar3 $ and $\bar4$.  The tangle
for $g_2$ has $2$ and $\bar5$ joined overcrossing a strand from $5$ and
$\bar2$ with three more vertical strands from the remaining nodes.  With this
it is straightforward to compute $e_4e_3g_2e_4$ as the tangle with $4$ and $2$
joined as well as $\bar 4$ and $\bar 2$.  Also $1$ and $\bar 1$ are joined
with a vertical line.  There are two more strands joining $5$ with $\bar3$ and
$3$ with $\bar5$ with the first strand overcrossing the second.  Now the
tangle $g_1e_4e_3g_2e_4g_1$ can be easily computed as the tangle with $4$ and
$2$ connected as well as $\bar4$ and $\bar2$ directly.  There are three
remaining strands which do not intersect these.  The first goes from $5$ to
$\bar1$.  The next goes from $3$ to $\bar3$ and passes under the first strand
crossing once.  The last strand connects $1$ with $\bar5$ and passes under
these two strands with two crossings.  The tangle for
$e_4e_3g_2e_4g_1e_4e_3g_2e_4$ is the same except there is an internal cycle
connecting $4,2,\bar2,\bar4$.  This gives the $\delta$ mentioned.  It is
straightforward to check that this relation handles all of the cases involving
$\hat{s_0}$ and $\hat{s_1}$ by showing $\hat{s}_0\hat{s_1}\hat{s}_0\isog
\hat{s}_1\hat{s_0}\hat{s}_1$.  For example
$g_1e_6e_5e_4e_3g_2e_4e_5e_6g_1\isog e_6e_5g_1e_4e_3g_2e_4g_1e_5e_6$.  Now use
$g_1e_4e_3g_2e_4g_1\isog \delta e_4e_3g_2e_4g_1e_4e_3g_2e_4$.

\end{proof}

\begin{Def}\label{def:HY}
\rm For each $Y$ of Table~\ref{table:types2} column~$5$ we let $\He_Y$ be the Hecke algebra generated by $\hat{s}_i$ mod $J_{t+1}$ as in
Proposition~\ref{hecke}.  Here $t$ is the size of the admissible
closure of $Y$ listed in the second column.
\end{Def}

We now describe the cell datum.  Fix $n\in\{6,7, 8\}$ and consider $M =
\E_n$.  For each $Y$ in Table~\ref{table:types2} column~$5$ for that $n$, we
let $(\Lambda_Y,D_Y,C_Y,*_Y)$ be the cell datum for the Hecke algebra
$\He_{Y}$ of type $M_{Y}$ listed in the fourth column for $Y$ as given by
Definition~\ref{def:HY}.  Here $t$ is the size of $B_Y$ listed in the second
column.  For $x,y \in D_Y$, $C_Y(x,y)$ is a coset mod $J_{t+1}$.  We would
like to have elements of $\BMW(M)$.  Each is a linear combination of words
in $\hat{s_i}$ and we can take the words in $J_t$ and not in $J_{t+1}$ if we
wish.  We define $C(x,y)$ as this sum.

Taken mod $J_{t+1}$ they are in $\He_Y$.

{From} \cite{MG} we know we can take $*_Y$ to be $\cdot ^\op$ for the Hecke
algebra.  Here, we let $*_Y$ be the restriction to $\He_{Y}$ of ${\cdot}^\op$
acting on the inverse image of $\He_Y$ in $\BMW(M)$.  Note that
${\cdot}^\op$ acts on $J_{t+1}$ and so acts on $\He_Y$. By \cite{MG}, these
cell data are known to exist if $S$ has inverses of the bad primes.  We
take the values of $C_Y$ in $\BMW(M)$ for each $Y\in \Lambda$ as discussed
above.  We want one more Hecke algebra for $Y=\emptyset$ which does not appear
in Table~\ref{table:types2}.  Here the Hecke algebra is $\BMW(M)$ mod
$J_1$.  Indeed $\BMW(M)/J_1$ is the Hecke algebra of type $M$.  We
denote this $\He_\emptyset$.  The braid relations are satisfied by definition
and the quadratic relations hold by (RSrr) as $e_i\in J_1$.  We let
$\Lambda_\emptyset$ be the poset for the cell datum for this Hecke algebra of
type $M$.  It it were in the table it would have $|X|=t=0$.

The poset $\Lambda$ is the disjoint union of the posets $\Lambda_Y$ of the
cell data for the various Hecke algebras $\He_Y$ together with $\Lambda_\emptyset$
for  $Y=\emptyset$.    We make $\Lambda$ into a poset as follows.
For a fixed $Y$, $\Lambda_Y$ it is already a poset, and we keep the same
partial order.  Furthermore, any element of $\Lambda_Y$ is greater than any
element of $\Lambda_{Y'}$ if $t<t'$ where $t'$ is the integer in column
two for the row of $Y'$.  This is the size of the admissible closure
of $Y'$.  In particular the elements of $\Lambda_\emptyset$
are greater than the elements of $\Lambda_Y$ for any $Y\neq \emptyset$.

 For
$\lambda \in \Lambda_Y$, we set $D(\lambda)= WX \times D_Y(\lambda)$ where $X$ is the admissible closure
of $Y$ whose size is listed in column~$2$ of Table~\ref{table:types2}. This determines $D$.  We identify $D(\Lambda_\emptyset)$ as just $\He_\emptyset$.

For a fixed $Y$ recall we have defined elements $a_{B}$ in
Definition~\ref{def:aB}. To distinguish the various choices of $Y$ we let
$a_{B,Y}$ be this element.  We now define words $\hat{a}_{B,Y}$ as the
natural elements of $\BMW(M)$.

\begin{Def} \rm
For each $B$ we make a choice of one of the words $a_{B,Y}$ given in Definition~\ref{def:aB}.  We then let
$\hat{a}_{B,Y}$ be $\rho(a_{B,Y})$.
\end{Def}

\np
We define $C$ as follows.  For
$\lambda\in\Lambda_Y$,  and  $(B,x),(B',y)\in D(\lambda)$, we have
$$C\big((B,x),(B',y)\big)=\hat{a}_{B,Y}C_Y(x,y)\hat{a}_{B',Y}^\op.$$

Since we already defined $*$ by the opposition map,
this concludes the definition of $(\Lambda, D, C, *)$.

\begin{Thm}\label{th:cellular}
Let $M$ be a spherical simply laced Coxeter type.  Let $S$ be an integral
domain containing $R$ with $p^{-1}\in S$ whenever $p$ is a bad prime for
$M$.  Then the quadruple $(\Lambda, T, C, *)$ is a cell datum for
$\BMW(M)\otimes_R S$, and so this algebra is cellular.
\end{Thm}

\begin{proof} Cellularity is known for $M= \ddA_n$ $(n\ge1)$ by \cite{Xi}
and for $M = \ddD_n$ $(n\ge4)$ by \cite{CGW3}.  By standard arguments it
remains to verify the conditions (C1), (C2), (C3) for
$M\in\{\ddE_6,\ddE_7,\ddE_8\}$.

\nl(C1) The map $C$ has been chosen so that its image is the set of all
$\hat{a}_{B,Y}C_Y(x,y)\hat{a}_{B',Y}$ where $Y\in \MY$ and $C_Y(x,y)$ are
elements of a basis of the Hecke algebra $\He_Y$.  This is a spanning set.
Injectivity follows from the ranks of the various quotients.

\nl(C2).
Clearly, $*=\cdot^\op$ is an $S$-linear anti-involution.
For each $Y$, choose
$\lambda\in\Lambda_Y$, and $(B,x),(B',y)\in D(\lambda)$.
Then $(\hat{a}_{B,Y}C_Y(x,y)\hat{a}_{B',Y}^\op)^\op
=\hat{a}_{B',Y}C_Y(x,y)^\op \hat{a}_{B,Y}^\op$,  so, in order to establish
$\big(C((B,x),(B',y))\big)^*=C((B',y),(B,x))$,
it suffices to verify
that $C_Y(x,y)^\op $ coincides with $C_Y(y,x)$.
Now $*_Y$ on
$\He_Y{(Y)}$  coincides with opposition, so
modulo $J_{t+1}$ we have $C_Y(x,y)^\op =  C_Y(x,y)^{*_Y} = C_Y(y,x)$
by the cellularity of $(\Lambda_Y,D_Y,C_Y,*_Y)$.
On the other hand, as the inverse image in $\BMW(M)\otimes_R S$ of
$\He_{Y}$ is invariant under opposition, and contains the values of $C_Y$,
it contains $C_Y(x,y)^\op - C_Y(y,x)$, so
$C_Y(x,y)^\op-C_Y(y,x)\in  J_{t+1} .$  However the elements of $C_Y$ were chosen
in $J_t\setminus J_{t+1}$ and so
$C_Y(x,y)^\op=C_Y(y,x)$, as required.

\nl (C3).  Let $\lambda\in\Lambda_Y$ and $(B,x),(B',y)\in D(\lambda)$.
Fix $Y$.
It clearly suffices to prove the formulas for $a$ running over the generators
$g_i$ and $e_i$ of $\BMW(M)\otimes_R S$.

By choice of $C_Y$, we have
$C_Y(x,y)$ in contained in the ideal generated by $\rho(e_{Y})$.  Using Theorem~\ref{th:aBcharacterization}, there is $h_{B,i}\in H_{Y}$,
depending only on $B$ and $i$,
such that $g_i \hat{a}_{B,Y} \in \hat{a}_{r_iB,Y}\rho(h_{B,i}) + J_{t+1}$ .
As $(\Lambda_Y,D_Y,C_Y,*_Y)$ is a cell datum for
$\He_{Y}$ mod $J_{t+1}$,
there are $\nu_{i}(u,B,x)\in S$, independent of $B'$ and $y$,
for each $u\in D_Y(\lambda)$
such that
\begin{eqnarray*}
\rho(h_{B,Y})C_Y(x,y)&\in&\sum_{u \in D_Y(\lambda)}
\nu_{i}(u,B,x)C_Y(u,y) + (\He_{Y})_{<\lambda} + J_{t+1}.
\end{eqnarray*}
Since both  $(\He_{Y})_{<\lambda}$  and $J_{t+1}$ are contained in
$\alg_{<\lambda}$, we find
\begin{eqnarray*}
g_iC((B,x),(B',y)) & = &
g_i\hat{a}_{B,Y}C_Y(x,y)\hat{a}_{B',Y}^\op\\
&\in&\hat{a}_{r_iB,Y}\rho(h_{B,i})C_Y(x,y)\hat{a}_{B',Y}^\op + \alg_{<\lambda}\\
&=&\sum_{u\in D_Y(\lambda)}
\nu_i(u,B,x)\hat{a}_{r_iB,Y}C_Y(u,y)\hat{a}_{B',Y}^\op + \alg_{<\lambda}\\
&=&\sum_{u\in D_Y(\lambda)}
\nu_i(u,B,x)C((r_iB,u),(B',y)) + \alg_{<\lambda}
\end{eqnarray*}
as required.

Rewriting (RSrr) to $e_i = lm^{-1}(g_i^2+mg_i-1)$, we see that, if $m^{-1}\in
S$, the proper behavior of the cell data under left multiplication by
$e_i$ is taken care of by the above formulae for $g_i$.  Otherwise
a proof using $e_i$ works just as above for $g_i$ again using Theorem~\ref{th:aBcharacterization}.

This establishes that $(\Lambda,T,C,*)$ is a cell datum for $\alg$
and so completes the proof of cellularity of $\BMW(M)\otimes_R S$.
\end{proof}

\begin{Remark}\label{substructures}
\rm Let $K$ be any set of nodes of $M$.  A consequence of Theorem
\ref{th:main} is that the standard parabolic subalgebra of type $K$, that is,
the subalgebra generated by $\{g_i,e_i\mid i\in K\}$
is naturally isomorphic
to the BMW algebra whose type is the restriction of $M$ to $K$.
\end{Remark}

\end{section}

\end{document}